\documentclass[12pt]{article}

\usepackage{amssymb}
\usepackage{bbm}
\usepackage{eufrak}

\newcommand{\R}{\mathbb{R}}

\newcommand{\Z}{\mathbb{Z}}
\renewcommand{\P}{\mathbb{P}}

\newcommand{\E}{\mathbb{E}}

\newcommand{\D}{\mathbb{D}}
\newcommand{\N}{\mathbb{N}}

\newcommand{\T}{\mathbb{T}}
\newcommand{\A}{\mathbb{A}}
\renewcommand{\Z}{\mathbb{Z}}
\renewcommand{\L}{\mathbb{L}}

\def\build#1_#2^#3{\mathrel{
\mathop{\kern 0pt#1}\limits_{#2}^{#3}}}

\def\cq{$\hfill \square$}

\def\ind{{\mathbbm{1}}_}

\def\al{\alpha}
\def\vep{\varepsilon}

\def\vph{\varphi}
\def\lm{\lambda}
\def\sg{\sigma}
\def\dl{\delta}

\def\l{{\cal L}}

\def\a{{\cal A}}
\def\m{{\cal M}}
\def\n{{\cal N}}
\def\g{{\cal G}}
\def\f{{\cal F}}

\def\h{{\cal H}}
\def\t{{\cal T}}

\def\x{{\cal X}}
\def\z{{\cal Z}}
\def\y{{\cal Y}}

\def\proj{{\mathfrak p}}
\def\tt{{\Theta}}
\def\PP{\hbox{\bf P}}

\def\EE{\hbox{\bf E}}

\def\XX{\hbox{\bf X}}

\def\TT{\mbox{\boldmath$\mathcal{T}$}}
\def\OO{\mbox{\boldmath$\Omega$}}

\def\ttt{\mbox{\boldmath$\theta$}}

\def\be{\begin{equation}}
\def\ee{\end{equation}}
\def\ba{\begin{eqnarray*}}
\def\ea{\end{eqnarray*}}

\def\wh{\widehat}
\def\wt{\widetilde}

\def\la{\longrightarrow}
\def\lmt{\longmapsto}
\def\da{\downarrow}
\def\ua{\uparrow}

\def\proof{\noindent{\textsc {Proof}\hspace{0.05cm}:\hspace{0.15cm}}}

\newtheorem{theorem}{Theorem}[section]
\newtheorem{lemma}[theorem]{Lemma}
\newtheorem{proposition}[theorem]{Proposition}
\newtheorem{corollary}[theorem]{Corollary}

\newtheorem{definition}{Definition}[section]

\begin{document}

\title{ \bf Regenerative real trees}
\author{Mathilde {\sc Weill} \\
{\small D.M.A., \'Ecole normale sup\'erieure, 45 rue d'Ulm, 75005 Paris, France}}
\vspace{2mm}
\date{\today}

\maketitle

\begin{abstract}
In this work, we give a description of all $\sg$-finite measures on the space of rooted compact $\R$-trees which satisfy a certain regenerative property. We show that any infinite measure which satisfies the regenerative property is the "law" of a L\'evy tree, that is, the "law" of a tree-valued random variable that describes the genealogy of a population evolving according to a continuous-state branching process. On the other hand, we prove that a probability measure with the regenerative property must be the law of the genealogical tree associated with a continuous-time discrete-state branching process.  
\end{abstract}

\section{Introduction}

Galton-Watson trees are well known to be characterized among all random discrete trees by a regenerative property. More precisely, if $\gamma$ is a probability measure on $\Z_+$, the law $\Pi_\gamma$ of the Galton-Watson tree with offspring distribution $\gamma$ is uniquely determined by the following two conditions: Under the probability measure $\Pi_\gamma$,
\begin{description}
\item (i) the ancestor has $p$ children with probability $\gamma(p)$,
\item (ii) if $\gamma(p)>0$, then conditionally on the event that the ancestor has $p$ children, the $p$ subtrees which describe the genealogy of the descendants of these children, are independent and distributed according to $\Pi_\gamma$.
\end{description} 

The aim of this work is to study $\sg$-finite measures satisfying an analogue of this property on the space of equivalence classes of rooted compact $\R$-trees.

An $\R$-tree is a metric space $(\t,d)$ such that for any two points $\sg_1$ and $\sg_2$ in $\t$, there is a unique arc with endpoints $\sg_1$ and $\sg_2$, and furthermore this arc is isometric to a compact interval of the real line. In this work, all $\R$-trees are supposed to be compact. A rooted $\R$-tree is an $\R$-tree with a distinguished vertex called the root. Say that two rooted $\R$-trees are equivalent if there is a root-preserving isometry that maps one onto the other. It was noted in \cite{EPW} that the set $\T$ of all equivalence classes of rooted compact $\R$-trees equipped with the pointed Gromov-Hausdorff distance $d_{GH}$ (see e.g. Chapter 7 in \cite{BBI}), is a Polish space. Hence it is legitimate to consider random variables with values in $\T$, that is, random $\R$-trees. A particularly important example is the CRT, which was introduced by Aldous \cite{Al1}, \cite{Al3} with a different formalism. Striking applications of the concept of random $\R$-trees can be found in the recent papers \cite{EPW} and \cite{EW}.

Let $\t$ be an $\R$-tree. We write $\h(\t)$ for the height of the $\R$-tree $\t$, that is the maximal distance from the root to a vertex of $\t$. For every $t\geq0$, we denote by $\t_{\leq t}$ the set of all vertices of $\t$ which are at distance at most $t$ from the root, and by $\t_{>t}$ the set of all vertices which are at distance greater than $t$ from the root. To each connected component of $\t_{>t}$ there corresponds a "subtree" of $\t$ above level $t$ (see section \ref{sousarbres} for a more precise definition). For every $h>0$, we define $Z(t,t+h)(\t)$ as the number of subtrees of $\t$ above level $t$ with height greater than $h$.   

Let $\tt$ be a $\sg$-finite measure on $\T$, such that $0<\tt(\h(\t)>t)<\infty$ for every $t>0$ and $\tt(\h(\t)=0)=0$. We say that $\tt$ satisfies the regenerative property (R) if the following holds:
\begin{description}
\item (R) For every $t,h>0$ and $p\in\N$, under the probability measure $\tt(\cdot\mid\h(\t)>t)$ and conditionally on the event $\{Z(t,t+h)=p\}$, the $p$ subtrees of $\t$ above level $t$ with height greater than $h$ are independent and distributed according to the probability measure $\tt(\cdot\mid\h(\t)>h)$.
\end{description}
This is a natural analogue of the regenerative property stated above for Galton-Watson trees. Beware that, unlike the discrete case, there is no natural order on the subtrees above a given level. So, the preceding property should be understood in the sense that the unordered collection of the $p$ subtrees in consideration is distributed as the unordered collection of $p$ independent copies of $\tt(\cdot\mid\h(\t)>h)$.

Property (R) is known to be satisfied by a wide class of infinite measures on $\T$, namely the "laws" of L\'evy trees. L\'evy trees have been introduced by T. Duquesne and J.F. Le Gall in \cite{DuLG}. Their precise definition is recalled in section \ref{secLevytree}, but let us immediately give an informal presentation.

Let $Y$ be a critical or subcritical continuous-state branching process.  The distribution of $Y$ is characterized by its branching mechanism function $\psi$. Assume that $Y$ becomes extinct a.s., which is equivalent to the condition $\int_1^\infty\psi(u)^{-1}du<\infty$. The $\psi$-L\'evy tree is a random variable taking values in $(\T,d_{GH})$, which describes the genealogy of a population evolving according to $Y$ and starting with infinitesimally small mass. More precisely, the "law" of the L\'evy tree is defined in \cite{DuLG} as a $\sg$-finite measure on the space $(\T,d_{GH})$, such that $0<\tt_\psi(\h(\t)>t)<\infty$ for every $t>0$. As a consequence of Theorem 4.2 of \cite{DuLG}, the measure $\tt_\psi$ satisfies Property (R). In the special case $\psi(u)=u^\al$, $1<\al\leq2$ corresponding to the so-called stable trees, this was used by Miermont \cite{M1}, \cite{M2} to introduce and to study certain fragmentation processes.

In the present work we describe all $\sg$-finite measures on $\T$ that satisfy Property (R). We show that the only infinite measures satisfying Property (R) are the measures $\tt_\psi$ associated with L\'evy trees. On the other hand, if $\tt$ is a finite measure satisfying Property (R), we can obviously restrict our attention to the case $\tt(\T)=1$ and we obtain that $\tt$ must be the law of the genealogical tree associated with a continuous-time discrete-state branching process.  

\begin{theorem}\label{casinfini}
Let $\tt$ be an infinite measure on $(\T,d_{GH})$ such that $\tt(\h(\t)=0)=0$ and $0<\tt(\h(\t)>t)<+\infty$ for every $t>0$. Assume that $\tt$ satisfies Property (R). Then, there exists a continuous-state branching process, whose  branching mechanism is denoted by $\psi$, which becomes extinct almost surely, such that $\tt=\tt_\psi$.
\end{theorem}

\begin{theorem}\label{casfini}
Let $\tt$ be a probability measure on $(\T,d_{GH})$ such that $\tt(\h(\t)=0)=0$ and $0<\tt(\h(\t)>t)$ for every $t>0$. Assume that $\tt$ satisfies Property (R). Then there exists $a>0$ and a critical or subcritical probability measure $\gamma$ on $\Z_+\backslash\{1\}$ such that $\tt$ is the law of the genealogical tree for a discrete-space continuous-time branching process with offspring distribution $\gamma$, where branchings occur at rate $a$.
\end{theorem} 
In other words, $\tt$ in Theorem \ref{casfini} can be described in the following way: There exists a real random variable $J$ such that under $\tt$: 
\begin{description}
\item (i) $J$ is distributed according to the exponential distribution with parameter $a$ and there exists $\sg_J\in\t$ such that $\t_{\leq J}=[[\rho,\sg_J]]$,
\item (ii) the number of subtrees above level $J$ is distributed according to $\gamma$ and is independent of $J$,
\item (iii) for every $p\geq2$, conditionally on $J$ and given the event that the number of subtrees above level $J$ is equal to $p$, these $p$ subtrees are independent and distributed according to $\tt$.   
\end{description} 

Theorem \ref{casinfini} is proved in section \ref{sec3}, after some preliminary results have been established in section \ref{prelim}. A key idea of the proof is to use the regenerative property to embed discrete Galton-Watson trees in our random real trees (Lemma \ref{arbreGW}). A technical difficulty comes from the fact that real trees are not ordered whereas Galton-Watson trees are usually defined as random ordered discrete trees (cf subsection \ref{secarbres} below). To overcome this difficulty, we assign a random ordering to the discrete trees embedded in real trees. Another major ingredient of the proof of Theorem \ref{casinfini} is the construction of a "local time" $L_t$ at every level $t$ of a random real tree governed by $\tt$. The local time process is then shown to be a continuous-state branching process with branching mechanism $\psi$, which makes it possible to identify $\tt$ with $\tt_\psi$. Theorem \ref{casfini} is proved in section \ref{sec4}. Several arguments are similar to the proof of Theorem \ref{casinfini}, so that we have skipped some details.

\section{Preliminaries}\label{prelim}

In this section, we recall some basic facts about branching processes, $\R$-trees and L\'evy trees. 

\subsection{Branching processes}

\subsubsection{Continuous-state branching processes}

A (continuous-time) continuous-state branching process (in short a CSBP) is a Markov process $Y=(Y_t,t\geq0)$ with values in the positive half-line $[0,+\infty)$, with a Feller semigroup $(Q_t,t\geq0)$ satisfying the following branching property: For every $t\geq0$ and $x,x'\geq0$,
$$Q_t(x,\cdot)\ast Q_t(x',\cdot)=Q_t(x+x',\cdot).$$
Informally, this means that the union of two independent populations started respectively at $x$ and $x'$ will evolve like a single population started at $x+x'$.

We will consider only the critical or subcritical case, meaning that, for every $t\geq0$ and $x\geq0$,
$$\int_{[0,+\infty)}yQ_t(x,dy)\leq1.$$
Then, if we exclude the trivial case where $Q_t(x,\cdot)=\dl_0$ for every $t>0$ and $x\geq0$, the Laplace functional of the semigroup can be written in the following form: For every $\lm\geq0$,
$$\int_{[0,+\infty)}e^{-\lm y}Q_t(x,dy)=\exp(-xu(t,\lm)),$$
where the function $(u(t,\lm),t\geq0,\lm\geq0)$ is determined by the differential equation
$$\frac{du(t,\lm)}{dt}=-\psi(u(t,\lm)),\;\;u(0,\lm)=\lm,$$
and $\psi:\R_+\la\R_+$ is of the form
\be\label{psi}
\psi(u)=\al u+\beta u^2+\int_{(0,+\infty)}(e^{-ur}-1+ur)\pi(dr),
\ee
where $\al,\beta\geq0$ and $\pi$ is a $\sg$-finite measure on $(0,+\infty)$ such that $\int_{(0,+\infty)}(r\wedge r^2)\pi(dr)<\infty$. The process $Y$ is called the $\psi$-continuous-state branching process (in short the $\psi$-CSBP). 

Continuous-state branching processes may also be obtained as weak limits of rescaled Galton-Watson processes. We recall that an offspring distribution is a probability measure on $\Z_+$. An offspring distribution $\mu$ is said to be critical if $\sum_{i\geq0}i\mu(i)=1$ and subcritical if $\sum_{i\geq0}i\mu(i)<1$. Let us state a result that can be derived from \cite{Gr} and \cite{La}.

\begin{theorem}\label{Grimvall}
Let $(\mu_n)_{n\geq1}$ be a sequence of offspring distributions. For every $n\geq1$, denote by $X^n$ a Galton-Watson process with offspring distribution $\mu_n$, started at $X^n_0=n$. Let $(m_n)_{n\geq1}$ be a nondecreasing sequence of positive integers converging to infinity. We define a sequence of processes $(Y^n)_{n\geq1}$ by setting, for every $t\geq0$ and $n\geq1$,
$$Y^n_t=n^{-1}X^n_{[m_nt]}.$$
Assume that, for every $t\geq0$, the sequence $(Y^n_t)_{n\geq1}$ converges in distribution to $Y_t$ where $Y=(Y_t,t\geq0)$ is an almost surely finite process such that $\P(Y_\dl>0)>0$ for some $\dl>0$. Then, $Y$ is a continuous-state branching process and the sequence of processes $(Y^n)_{n\geq1}$ converges to $Y$ in distribution in the Skorokhod space $\D(\R_+)$.
\end{theorem}

\proof It follows from the proof of Theorem 1 of \cite{La} that $Y$ is a CSBP. Then, thanks to Theorem 2 of \cite{La}, there exists a sequence of offspring distributions $(\nu_{n})_{n\geq1}$ and a nondecreasing sequence of positive integers $(c_n)_{n\geq1}$ such that we can construct for every $n\geq1$ a Galton-Watson process $Z^n$ started at $c_n$ and with offspring distribution $\nu_{n}$ satisfying
$$\left(c_n^{-1}Z^n_{[nt]},t\geq0\right)\build{\la}_{n\to\infty}^{d}(Y_t,t\geq0),$$
where the symbol $\build{\la}_{}^{d}$ indicates convergence in distribution in $\D(\R_+)$. 

Let $(m_{n_k})_{k\geq1}$ be a strictly increasing subsequence of $(m_n)_{n\geq1}$. For $n\geq1$, we set $B^n=X^{n_k}$ and $b_n=n_k$ if $n=m_{n_k}$ for some $k\geq1$, and we set $B^n=Z^n$ and $b_n=c_n$ if there is no $k\geq1$ such that $n=m_{n_k}$. Then, for every $t\geq0$, $(b_n^{-1}B^n_{[nt]})_{n\geq1}$ converges in distribution to $Y_t$. Applying Theorem 4.1 of \cite{Gr}, we obtain that
$$\left(b_n^{-1}B^n_{[nt]},t\geq0\right)\build{\la}_{n\to\infty}^{d}(Y_t,t\geq0).$$
In particular, we have,
\be\label{soussuite}
\left(Y^{n_k}_t,t\geq0\right)\build{\la}_{k\to\infty}^{d}(Y_t,t\geq0).
\ee
As (\ref{soussuite}) holds for every strictly increasing subsequence of $(m_n)_{n\geq1}$, we obtain the desired result.\cq 

\subsubsection{Discrete-state branching processes}

A (continuous-time) discrete-state branching process (in short DSBP) is a continuous-time Markov chain $Y=(Y_t,t\geq0)$ with values in $\Z_+$ whose transition probabilities $(P_t(i,j),t\geq0)_{i\geq0,j\geq0}$ satisfy the following branching property: For every $i\in\Z_+$, $t\geq0$ and $|s|\leq1$,
$$\sum_{j=0}^\infty P_t(i,j)s^j=\left(\sum_{j=0}^\infty P_t(1,j)s^j\right)^i.$$
We exclude the trivial case where $P_t(i,i)=1$ for every $t\geq0$ and $i\in\Z_+$. Then, there exist $a>0$ and an offspring distribution $\gamma$ with $\gamma(1)=0$ such that the generator of $Y$ can be written of the form
$$Q=\left(\begin{array}{cccccc} 0 & 0 & 0 & 0 & 0 & \ldots \\
                               a\gamma(0) & -a & a\gamma(2) & a\gamma(3) & a\gamma(4) & \ldots \\
                               0 & 2a\gamma(0) & -2a & 2a\gamma(2) & a\gamma(3) & \ldots \\
                               0 & 0 & 3a\gamma(0) & -3a& 3a\gamma(2) & \ldots \\
                               \vdots & \vdots & \ddots & \ddots & \ddots & \ddots
                               \end{array}\right).$$
Furthermore, it is well known that $Y$ becomes extinct almost surely if and only if $\gamma$ is critical or subcritical. We refer the reader to \cite{AN} and \cite{No} for more details.

\subsection{Deterministic trees}
\subsubsection{The space $(\T,d_{GH})$ of rooted compact $\R$-trees}\label{arbresreels} 

We start with a basic definition.
\begin{definition}
A metric space $(\t,d)$ is an $\R$-tree if the following two properties hold for every $\sg_1,\sg_2\in\t$.
\begin{description}
\item (i) There is a unique isometric map $f_{\sg_1,\sg_2}$ from $[0,d(\sg_1,\sg_2)]$ into $\t$ such that 
$$f_{\sg_1,\sg_2}(0)=\sg_1\;\;and\;\; f_{\sg_1,\sg_2}(d(\sg_1,\sg_2))=\sg_2.$$
\item (ii) If $q$ is a continuous injective map from $[0,1]$ into $\t$ such that $q(0)=\sg_1$ and $q(1)=\sg_2$, we have
$$q([0,1])=f_{\sg_1,\sg_2}([0,d(\sg_1,\sg_2)]).$$
\end{description}
A rooted $\R$-tree is an $\R$-tree with a distinguished vertex $\rho=\rho(\t)$ called the root. 
\end{definition}
In what follows, $\R$-trees will always be rooted.

Let $(\t,d)$ be an $\R$-tree with root $\rho$, and $\sg,\sg_1,\sg_2\in\t$. We write $[[\sg_1,\sg_2]]$ for the range of the map $f_{\sg_1,\sg_2}$. In particular, $[[\rho,\sg]]$ is the path going from the root to $\sg$ and can be interpreted as the ancestral line of $\sg$. 

The height $\h(\t)$ of the $\R$-tree $\t$ is defined by $\h(\t
)=\sup\{d(\rho,\sg):\sg\in\t\}$. In particular, if $\t$ is compact, its height $\h(\t)$ is finite.

Two rooted $\R$-trees $\t$ and $\t'$ are called equivalent if there is a root-preserving isometry that maps $\t$ onto $\t'$. We denote by $\T$ the set of all equivalence classes of rooted compact $\R$-trees. We often abuse notation and identify a rooted compact $\R$-tree with its equivalence class. 

The set $\T$ can be equipped with the pointed Gromov-Hausdorff distance, which is defined as follows. If $(E,\dl)$ is a metric space, we use the notation $\dl_{Haus}$ for the usual Hausdorff metric between compact subsets of $E$. Then, if $\t$ and $\t'$ are two rooted compact $\R$-trees with respective roots $\rho$ and $\rho'$, we define the distance $d_{GH}(\t,\t')$ as
$$d_{GH}(\t,\t')=\inf\Big\{\dl_{Haus}\left(\phi(\t),\phi'(\t')\right)\vee\dl(\phi(\rho),\phi'(\rho'))\Big\},$$
where the infimum is over all isometric embeddings $\phi:\t\la E$ and $\phi':\t'\la E$ into a common metric space $(E,\dl)$. We see that $d_{GH}(\t,\t')$ only depends on the equivalence classes of $\t$ and $\t'$. Furthermore, according to Theorem 2 in \cite{EPW}, $d_{GH}$ defines a metric on $\T$ that makes it complete and separable. Notice that $d_{GH}(\t,\t')$ makes sense more generally if $\t$ and $\t'$ are pointed compact metric spaces (see e.g. Chapter 7 in \cite{BBI}). We will use this in the proof of Lemma \ref{discretisation} below.

We equip $\T$ with its Borel $\sg$-field. If $\t\in\T$, we set $\t_{\leq t}=\{\sg\in\t:d(\rho,\sg)\leq t\}$ for every $t\geq0$. Plainly, $\t_{\leq t}$ is a real tree whose root is the same as the root of $\t$. Note that the mapping $\t\lmt\t_{\leq t}$ from $\T$ into $\T$ is Lipschitz for the Gromov-Hausdorff metric. 
  
\subsubsection{The $\R$-tree coded by a function}\label{arbrescodes}

We now recall a construction of rooted compact $\R$-trees which is described in \cite{DuLG}. Let $g:[0,+\infty)\la[0,+\infty)$ be a continuous function with compact support satisfying $g(0)=0$. We exclude the trivial case where $g$ is identically zero. For every $s,t\geq0$, we set
$$m_g(s,t)=\inf_{r\in[s\wedge t,s\vee t]}g(r)$$
and
$$d_g(s,t)=g(s)+g(t)-2m_g(s,t).$$
We define an equivalence relation $\sim$ on $[0,+\infty)$ by declaring that $s\sim t$ if and only if $d_g(s,t)=0$ (or equivalently if and only if $g(s)=g(t)=m_g(s,t)$). Let $\t_g$ be the quotient space 
$$\t_g=[0,+\infty)/\sim.$$
Then, $d_g$ induces a metric on $\t_g$ and we keep the notation $d_g$ for this metric. According to Theorem 2.1 of \cite{DuLG}, the metric space $(\t_g,d_g)$ is a compact $\R$-tree. By convention, its root is the equivalence class of $0$ for $\sim$ and is denoted by $\rho_g$.

\subsubsection{Subtrees of a tree above a fixed level}\label{sousarbres}

Let $(\t,d)\in\T$ and $t>0$. Denote by $\t^{i,\circ},i\in I$ the connected components of the open set $\t_{>t}=\{\sg\in\t:d(\rho(\t),\sg)>t\}$. Let $i\in I$. Then the ancestor of $\sg$ at level $t$, that is, the unique vertex on the line segment $[[\rho,\sg]]$ at distance $t$ from $\rho$, must be the same for all $\sg\in\t^{i,\circ}$ . We denote by $\sg_i$ this common ancestor and set $\t^{i}=\t^{i,\circ}\cup\{\sg_i\}$. Then $\t^{i}$ is a compact rooted $\R$-tree with root $\sg_i$. The trees $\t^{i},i\in I$ are called the subtrees of $\t$ above level $t$. We now consider, for every $h>0$,
$$Z(t,t+h)(\t)=\#\{i\in I:\h(\t^{i})>h\}.$$
By a compactness argument, we can easily verify that $Z(t,t+h)(\t)<\infty$.

\subsubsection{Discrete trees and real trees}\label{secarbres}

We start with some formalism for discrete trees. We first introduce the set of labels 
$$U=\bigcup_{n\geq0}\N^n,$$ 
where by convention $\N^0=\{\emptyset\}$. An element of $U$ is a sequence $u=u^1\ldots u^n$, and we set $|u|=n$ so that $|u|$ represents the generation of $u$. In particular, $\vert\emptyset\vert=0$. If $u=u^1\ldots u^n$ and $v=v^1\ldots v^m$ belong to $U$, we write $uv=u^1\ldots u^nv^1\ldots v^m$ for the concatenation of $u$ and $v$. In particular, $\emptyset u=u\emptyset=u$. The mapping $\pi:U\backslash\{\emptyset\}\la U$ is defined by $\pi(u^1\ldots u^n)=u^1\ldots u^{n-1}$ ($\pi(u)$ is the father of $u$). Note that $\pi^k(u)=\emptyset$ if $k=|u|$. 

A rooted ordered tree $\theta$ is a finite subset of $U$ such that 
\begin{description}
\item (i) $\emptyset\in\theta$,
\item (ii) $u\in\theta\backslash\{\emptyset\}\Longrightarrow\pi(u)\in\theta$,
\item (iii) for every $u\in\theta$, there exists a number $k_u(\theta)\geq0$ such that $uj\in\theta$ if and only if $1\leq j\leq k_u(\theta)$.
\end{description}
We denote by $\a$ the set of all rooted ordered trees. If $\theta\in \a$, we write $\h(\theta)$ for the height of $\theta$, that is $\h(\theta)=\max\lbrace\vert u\vert:u\in\theta\rbrace$. And for every $u\in\theta$, we define $\tau_u\theta\in\a$ by $\tau_u\theta=\{v\in U:uv\in\theta\}$. This is the tree $\theta$ shifted at $u$. 

Let us define an equivalence relation on $\a$ by setting $\theta\sim\theta'$ if and only if we can find a permutation $\vph_u$ of the set $\{1,\ldots,k_u(\theta)\}$ for every $u\in\theta$ such that $k_u(\theta)\geq1$, in such a way that
$$\theta'=\{\emptyset\}\cup\left\{\vph_\emptyset(u^1)\vph_{u^1}(u^2)\ldots\vph_{u^1\ldots u^{n-1}}(u^n):u^1\ldots u^n\in\theta,n\geq1\right\}.$$
In other words $\theta\sim\theta'$ if they correspond to the same unordered tree. Let $\A=\a/\sim$ be the associated quotient space and let $\proj:\a\la\A$ be the canonical projection. It is immediate that if $\theta\sim\theta'$, then $k_\emptyset(\theta)=k_\emptyset(\theta')$. So, for every $\xi\in\A$, we may define $k_\emptyset(\xi)=k_\emptyset(\theta)$ where $\theta$ is any representative of $\xi$. Let us fix $\xi\in\A$ such that $k_\emptyset(\xi)=k>0$ and choose a representative $\theta$ of $\xi$. We can define $\{\xi^1,\ldots,\xi^k\}=\{\proj(\tau_1\theta),\ldots,\proj(\tau_k\theta)\}$ as the unordered family of subtrees of $\xi$ above the first generation. Then, if $F:\a^k\la\R_+$ is any symmetric measurable function, we have
\begin{eqnarray}
&&\hspace{-1cm}\left(\#\proj^{-1}(\xi)\right)^{-1}\sum_{\theta\in\proj^{-1}(\xi)}F(\tau_1\theta,\ldots,\tau_k\theta)\nonumber\\
&=&\left(\#\proj^{-1}(\xi^1)\right)^{-1}\ldots\left(\#\proj^{-1}(\xi^k)\right)^{-1}\sum_{\theta_1\in\proj^{-1}(\xi^1)}\ldots\sum_{\theta_k\in\proj^{-1}(\xi^k)}F(\theta_1,\ldots,\theta_k).\label{card}
\end{eqnarray}
Note that the right-hand side of (\ref{card}) is well defined since it is symmetric in $\{\xi^1,\ldots,\xi^k\}$. The identity (\ref{card}) is a simple combinatorial fact, whose proof is left to the reader.

A marked tree is a pair $T=(\theta,\{h_u\}_{u\in\theta})$ where $\theta\in\a$ and $h_u\geq0$ for every $u\in\theta$. We denote by $\m$ the set of all marked trees. We can associate with every marked tree $T=(\theta,\{h_u\}_{u\in\theta})\in\m$, an $\R$-tree $\t^T$ in the following way. Let $\R^\theta$ be the vector space of all mappings from $\theta$ into $\R$. Write $(e_u,u\in\theta)$ for the canonical basis of $\R^\theta$. We define $l_\emptyset=0$ and $l_u=\sum_{k=1}^{|u|}h_{\pi^k(u)}e_{\pi^k(u)}$ for $u\in\theta$. Let us set
$$\t^T=\bigcup_{u\in\theta}[l_u,l_u+h_ue_u].$$
$\t^T$ is a connected union of line segments in $\R^\theta$. It is equipped with the distance $d_T$ such that $d_T(a,b)$ is the length of the shortest path in $\t^T$ between $a$ and $b$, and can be rooted in $\rho(\t^T)=0$ so that it becomes a rooted compact $\R$-tree.
 
If $\theta\in\a$, we write $\t^\theta$ for the $\R$-tree $\t^T$ where $T=(\theta,\{h_u\}_{u\in\theta})$ with $h_\emptyset=0$ and $h_u=1$ for every $u\in\theta\setminus\{\emptyset\}$, and we write $d_\theta$ for the associated distance. We then set $m_\emptyset=0$ and $m_u=\sum_{k=0}^{|u|-1}e_{\pi^k(u)}=l_u+e_u$ for every $u\in\theta\setminus\{\emptyset\}$. 

It is easily checked that $\t^\theta=\t^{\theta'}$ if $\theta\sim\theta'$. Thus for every $\xi\in\A$, we may write $\t^\xi$ for the tree $\t^\theta$ where $\theta$ is any representative of $\xi$. 

We will now explain how to approximate a general tree $\t$ in $\T$ by a discrete type tree. Let $\vep>0$ and set $\T^{(\vep)}=\{\t\in\T:\h(\t)>\vep\}$. For every $\t$ in $\T^{(\vep)}$, we can construct by induction an element $\xi^\vep(\t)$ of $\A$ in the following way:
\begin{description}
\item $\bullet$\hspace{0.5cm}If $\t\in\T^{(\vep)}$ satisfies $\h(\t)\leq2\vep$, we set $\xi^\vep(\t)=\proj(\{\emptyset\})$. 
\item $\bullet$\hspace{0.5cm}Let $n$ be a positive integer. Assume that we have defined $\xi^\vep(\t)$ for every $\t\in\T^{(\vep)}$ such that $\h(\t)\leq(n+1)\vep$. Let $\t$ be an $\R$-tree such that $(n+1)\vep<\h(\t)\leq(n+2)\vep$. We set $k=Z(\vep,2\vep)(\t)$ and we denote by $\t^1,\ldots,\t^k$ the $k$ subtrees of $\t$ above level $\vep$ with height greater than $\vep$. Then $\vep<\h(\t^i)\leq(n+1)\vep$ for every $i\in\{1,\ldots,k\}$, so we can define $\xi^\vep(\t^i)$. Let us choose a representative $\theta^i$ of $\xi^\vep(\t^i)$ for every $i\in\{1,\ldots,k\}$. We set
$$\xi^\vep(\t)=\proj\left(\{\emptyset\}\cup1\theta^1\cup\ldots\cup k\theta^k\right),$$
where $i\theta^i=\{iu:u\in\theta^i\}$. Clearly this does not depend on the choice of the representatives $\theta_i$.
\end{description}

If $r>0$ and $\t$ is a compact rooted $\R$-tree with metric $d$, we write $r\t$ for the same tree equipped with the metric $rd$.

\begin{lemma}\label{discretisation}
For every $\vep>0$ and every $\t\in\T^{(\vep)}$, we have
\be\label{eqdiscretisation}
d_{GH}(\vep\t^{\xi^\vep(\t)},\t)\leq4\vep. 
\ee
\end{lemma}

\proof Let $\vep>0$ and $\t\in\T$. Let $\theta$ be any representative of $\xi^\vep(\t)$. Recall the notation $(m_u,u\in\theta)$. We can construct a mapping $\phi:\theta\la\t$ such that:
\begin{description}
\item(i) For every $\sg\in\t$, there exists $u\in\theta$ such that $d(\sg,\phi(u))\leq2\vep$,
\item(ii) for every $u\in\theta$, $d(\phi(u),\rho)=\vep|u|$ where $\rho$ denotes the root of $\t$,
\item(iii) for every $u,u'\in\theta$, $0\leq \vep d_\theta(m_u,m_{u'})-d(\phi(u),\phi(u'))\leq2\vep$.
\end{description}
To be specific, we always take $\phi(\emptyset)=\rho$, which suffices for the construction if $\h(\t)\leq2\vep$. If $(n+1)\vep<\h(\t)\leq(n+2)\vep$ for some $n\geq1$, we have as above
$$\theta=\{\emptyset\}\cup1\theta^1\cup\ldots\cup k\theta^k,$$
where $\theta^1,\ldots,\theta^k$ are representatives of respectively $\xi^\vep(\t^1),\ldots,\xi^\vep(\t^k)$, if $\t^1,\ldots,\t^k$ are the subtrees of $\t$ above level $\vep$ with height greater than $\vep$. With an obvious notation we define $\phi(ju)=\phi_j(u)$ for every $j\in\{1,\ldots,k\}$ and $u\in\theta^j$. Properties (i)-(iii) are then easily checked by induction.

Let us now set $\Phi=\{\phi(u),u\in\theta\}$ and $M=\{m_u,u\in\theta\}$. 
We equip $\Phi$ with the metric induced by $d$ and $M$ with the metric induced by $\vep d_\theta$. Then, $\Phi$ and $M$ can be viewed as pointed compact metric spaces with respective roots $\rho$ and $0$. It is immediate that $d_{GH}(\vep\t^\theta,M)\leq\vep$. Furthermore, Property (i) above implies that $d_{GH}(\t,\Phi)\leq2\vep$. At last, according to Lemma 2.3 in \cite{EPW} and Property (iii) above, we have $d_{GH}(\Phi,M)\leq\vep$. Lemma \ref{discretisation} then follows from the triangle inequality for $d_{GH}$.\cq

\subsection{L\'evy trees}\label{secLevytree}

Roughly speaking, a L\'evy tree is a $\T$-valued random variable which is associated with a CSBP in such a way that it describes the genealogy of a population evolving according to this CSBP. 

\subsubsection{The measure $\tt_\psi$}

We consider on a probability space $(\OO,\PP)$ a $\psi$-CSBP $Y=(Y_t,t\geq0)$, where the function $\psi$ is of the form (\ref{psi}), and we suppose that $Y$ becomes extinct almost surely. This condition is equivalent to 
\be\label{condproclevy}
\int_1^\infty\frac{du}{\psi(u)}<\infty.
\ee 
This implies that at least one of the following two conditions holds:
\be\label{2cond}
\beta>0\;\;{\rm or}\;\;\int_{(0,1)}r\pi(dr)=\infty.
\ee

The L\'evy tree associated to $Y$ will be defined as the tree coded by the so-called height process, which is a functional of the L\'evy process with Laplace exponent $\psi$. Let us denote by $X=(X_t,t\geq0)$ a L\'evy process on $(\OO,\PP)$ with Laplace exponent $\psi$. This means that $X$ is a L\'evy process with no negative jumps, and that for every $\lm,t\geq0$,
$$\EE(\exp(-\lm X_t))=\exp(t\psi(\lm)).$$
Then, $X$ does not drift to $+\infty$ and has paths of infinite variation (by (\ref{2cond})). 

We can define the height process $H=(H_t,t\geq0)$ by the following approximation:
$$H_t=\lim_{\vep\to0}\frac{1}{\vep}\int_0^t\ind{\{X_s\leq I_t^s+\vep\}}ds,$$
where $I_t^s=\inf\{X_r:s\leq r\leq t\}$ and the convergence holds in probability (see Chapter 1 in \cite{DuLG-Mon}). Informally, we can say that $H$ measures the size of the set $\{s\in[0,t]:X_{s-}\leq I^s_t\}$. Thanks to condition (\ref{condproclevy}), we know that the process $H$ has a continuous modification (see Theorem 4.7 in \cite{DuLG-Mon}). From now on, we consider only this modification.

Let us now set $I_t=\inf{\{X_s:0\leq s\leq t\}}$ for every $t\geq0$, and consider the process $X-I=(X_t-I_t,t\geq0)$. We recall that $X-I$ is a strong Markov process, for which the point $0$ is regular. The process $-I$ is a local time for $X-I$ at level $0$. We write $N$ for the associated excursion measure. We let $\Delta(de)$ be the "law" of $(H_s,s\geq0)$ under $N$. This makes sense because the values of the height process in an excursion of $X-I$ away from $0$ only depend on that excursion (see section 1.2 in \cite{DuLG-Mon}). Then, $\Delta(de)$ is a $\sg$-finite measure on $C([0,\infty))$, and is supported on functions with compact support such that $e(0)=0$.

The L\'evy tree is the tree $(\t_e,d_e)$ coded by the function $e$, in the sense of section \ref{arbrescodes}, under the measure $\Delta(de)$. We denote by $\tt_\psi$ the $\sg$-finite measure on $\T$ which is the ``law''of the L\'evy tree, that is the image of $\Delta(de)$ under the mapping $e\lmt\t_e$.

\subsubsection{A discrete approximation of the L\'evy tree}\label{sectionGW}

Let us now recall that the L\'evy tree is the limit in the Gromov-Hausdorff distance of suitably rescaled Galton-Watson trees.

We start by recalling the definition of Galton-Watson trees which was given informally in the introduction above. Let $\gamma$ be a critical or subcritical offspring distribution. We exclude the trivial case where $\gamma(1)=1$. Then, there exists a unique probability measure $\Pi_\gamma$ on $\a$ such that:
\begin{description}
\item (i) For every $p\geq0$, $\Pi_\gamma(k_\emptyset=p)=\gamma(p)$,
\item (ii) for every $p\geq1$ with $\gamma(p)>0$, under the probability measure $\Pi_\gamma(\cdot\mid k_\emptyset=p)$, the shifted trees $\tau_1\theta,\ldots,\tau_p\theta$ are independent and distributed according to $\Pi_\gamma$. 
\end{description}

Recall that if $r>0$ and $\t$ is a compact rooted $\R$-tree with metric $d$, we write $r\t$ for the same tree equipped with the metric $rd$. The following result is Theorem 4.1 in \cite{DuLG}. 

\begin{theorem}\label{DuLG}
Let $(\gamma_n)_{n\geq1}$ be a sequence of critical or subcritical offspring distributions. For every $n\geq1$, let us denote by $X^n$ a Galton-Watson process with offspring distribution $\gamma_n$, started at $X^n_0=n$. Let $(m_n)_{n\geq1}$ be a nondecreasing sequence of positive integers converging to infinity. We define a sequence of processes $(Y^n)_{n\geq1}$ by setting, for every $t\geq0$ and $n\geq1$,
$$Y^n_t=n^{-1}X^n_{[m_nt]}.$$
Assume that, for every $t\geq0$, $(Y^n_t)_{n\geq1}$ converges in distribution to $Y_t$ where \linebreak $Y=(Y_t,t\geq0)$ is a $\psi$-CSBP which becomes extinct almost surely. Assume furthermore that for every $\dl>0$,
$$\liminf_{n\to\infty}\P(Y^n_\dl=0)>0.$$
Then, for every $a>0$, the law of the $\R$-tree $m_n^{-1}\t^\theta$ under $\Pi_{\gamma_n}(\cdot\mid\h(\theta)\geq[am_n])$ converges as $n\la\infty$ to the probability measure $\tt_\psi(\cdot\mid\h(\t)>a)$ in the sense of weak convergence of measures in the space $\T$. 
\end{theorem}

\section{Proof of Theorem \ref{casinfini}}\label{sec3}

Let $\tt$ be an infinite measure on $(\T,d_{GH})$ satisfying the assumptions of Theorem 1.1. Clearly $\tt$ is $\sg$-finite.

We start with two important lemmas that will be used throughout this section. Let us first define $v:(0,\infty)\la(0,\infty)$ by $v(t)=\tt(\h(\t)>t)$ for every $t>0$. For every $t>0$, we denote by $\tt^t$ the probability measure $\tt(\cdot\mid\h(\t)>t)$.  

\begin{lemma}\label{vcontinue}
The function $v$ is nonincreasing, continuous and verifies
$$ v(t)\build{\longrightarrow}_{t\to0}^{}\infty \;\;\;and\;\;\; v(t)\build{\longrightarrow}_{t\to\infty}^{}0.$$ 
\end{lemma}

\proof We only have to prove the continuity of $v$. To this end, we argue by contradiction and assume that there exists $t>0$ such that $\tt(\h(\t)=t)>0$. Let $s>0$ and $u\in(0,t)$ such that $v(u)>v(t)$. From the regenerative property (R), we have
\begin{eqnarray*}
\tt^s(\h(\t)=s+t)&=&\tt^s\left(\tt^s(\h(\t)=s+t\mid Z(s,s+u))\right)\\
&=&\tt^s\left(Z(s,s+u)\tt^u(\h(\t)=t)\left(\tt^u(\h(\t)\leq t)\right)^{Z(s,s+u)-1}\right)\\
&=&\frac{\tt(\h(\t)=t)}{v(u)}\;\tt^s\left(Z(s,s+u)\left(1-\frac{v(t)}{v(u)}\right)^{Z(s,s+u)-1}\right)\\
&>&0.
\end{eqnarray*}
We have shown that $\tt(\h(\t)=t+s)>0$ for every $s>0$. This is absurd since $\tt$ is $\sg$-finite.\cq

\begin{lemma}\label{binomiale}
For every $t>0$ and $0<a<b$, the conditional law of the random variable $Z(t,t+b)$, under the probability measure $\tt^t$ and given $Z(t,t+a)$, is a binomial distribution with parameters $Z(t,t+a)$ and $v(b)/v(a)$ (where we define the binomial distribution with parameters $0$ and $p\in[0,1]$ as the Dirac measure $\dl_0$).
\end{lemma}

\proof This is a straightforward consequence of the regenerative property.\cq

\subsection{The CSBP derived from $\tt$}

In this section, we consider a random forest of trees derived from a Poisson point measure with intensity $\tt$. We associate with this forest a family of Galton-Watson processes. We then construct local times at every level $a>0$ as limits of the rescaled Galton-Watson processes. Finally we show that the local time process is a CSBP.  
 
Let us now fix the framework. We consider a probability space $(\Omega,\P)$ and on this space a Poisson point measure $\n=\sum_{i\in I}\dl_{\t_i}$ on $\T$, whose intensity is the measure $\tt$. 

\subsubsection{A family of Galton-Watson trees}

We start with some notation that we need in the first lemma. We consider on another probability space $(\Omega',\P')$, a collection $(\theta_\xi,\xi\in\A)$ of independent $\a$-valued random variables such that for every $\xi\in\A$, $\theta_\xi$ is distributed uniformly over $\proj^{-1}(\xi)$. In what follows, to simplify notation, we identify an element $\xi$ of the set $\A$ with the subset $\proj^{-1}(\xi)$ of $\a$. Recall the definition of $\xi^\vep(\t)$ before Lemma \ref{discretisation}.

\begin{lemma}\label{arbreGW}
Let us define for every $\vep>0$, a mapping $\theta^{(\vep)}$ from $\T^{(\vep)}\times\Omega'$ into $\a$ by
$$\theta^{(\vep)}(\t,\omega)=\theta_{\xi^\vep(\t)}(\omega).$$
Then for every positive integer $p$, the law of the random variable $\theta^{(\vep)}$ under the probability measure $\tt^{p\vep}\otimes\P'$ is $\Pi_{\mu_\vep}(\cdot\mid\h(\theta)\geq p-1)$ where $\mu_\vep$ denotes the law of $Z(\vep,2\vep)$ under $\tt^\vep$.
\end{lemma} 

\proof Since $\{\h(\t)>p\vep\}\times\Omega'=\{\h(\theta^{(\vep)})\geq p-1\}$ for every $p\geq1$, it suffices to show the result for $p=1$. Let $k$ be a nonnegative integer. According to the construction of $\xi^\vep(\t)$, we have
$$\tt^\vep\otimes\P'\left(k_\emptyset\left(\theta^{(\vep)}\right)=k\right)=\tt^\vep(Z(\vep,2\vep)=k)=\mu_\vep(k).$$ 
Let us fix $k\geq1$ with $\mu_\vep(k)>0$. Let $F:\a^k\la\R_+$ be a symmetric measurable function. Then we have
\begin{eqnarray}
&&\hspace{-1cm}\tt^\vep\otimes\P'\left(F\left(\tau_1\theta^{(\vep)},\ldots,\tau_k\theta^{(\vep)}\right)\Big|\;k_\emptyset\left(\theta^{(\vep)}\right)=k\right)\nonumber\\
&=&\tt^\vep\otimes\P'\Big(\sum_{\theta\in\xi^\vep(\t)}F\left(\tau_1\theta,\ldots,\tau_k\theta\right)\ind{\{\theta_{\xi^\vep(\t)}=\theta\}}\Big|\;Z(\vep,2\vep)=k\Big)\nonumber\\
&=&\tt^\vep\Big((\#\xi^\vep(\t))^{-1}\sum_{\theta\in\xi^\vep(\t)}F\left(\tau_1\theta,\ldots,\tau_k\theta\right)\Big|\;Z(\vep,2\vep)=k\Big)\label{arbreGW3}.
\end{eqnarray}
On the event $\{Z(\vep,2\vep)=k\}$, we write $\t^1,\ldots,\t^k$ for the $k$ subtrees of $\t$ above level $\vep$ with height greater than $\vep$. Then, Formula (\ref{card}) and the regenerative property yield
\ba
&&\hspace{-1cm}\tt^\vep\Big((\#\xi^\vep(\t))^{-1}\sum_{\theta\in\xi^\vep(\t)}F\left(\tau_1\theta,\ldots,\tau_k\theta\right)\Big|\;Z(\vep,2\vep)=k\Big)\\
&=&\tt^\vep\Big((\#\xi^\vep(\t^1))^{-1}\ldots(\#\xi^\vep(\t^k))^{-1}\sum_{\theta_1\in\xi^\vep(\t^1)}\ldots\sum_{\theta_k\in\xi^\vep(\t^k)}F\left(\theta_1,\ldots,\theta_k\right)\Big|\;Z(\vep,2\vep)=k\Big)\\
&=&\int\tt^\vep(d\t_1)\ldots\tt^\vep(d\t_k)(\#\xi^\vep(\t_1))^{-1}\ldots(\#\xi^\vep(\t_k))^{-1}\sum_{\theta_1\in\xi^\vep(\t_1)}\ldots\sum_{\theta_k\in\xi^\vep(\t_k)}F\left(\theta_1,\ldots,\theta_k\right)\\
&=&\int\tt^\vep\otimes\P'(d\t_1,d\omega'_1)\ldots\tt^\vep\otimes\P'(d\t_k,d\omega'_k)F\left(\theta^{(\vep)}(\t_1,\omega'_1),\ldots,\theta^{(\vep)}(\t_k,\omega'_k)\right),
\ea
as in (\ref{arbreGW3}). We have thus proved that
\begin{eqnarray}
&&\hspace{-1cm}\tt^\vep\otimes\P'\left(F\left(\tau_1\theta^{(\vep)},\ldots,\tau_k\theta^{(\vep)}\right)\Big|\;k_\emptyset\left(\theta^{(\vep)}\right)=k\right)\nonumber\\
&=&\int\tt^\vep\otimes\P'(d\t_1,d\omega'_1)\ldots\tt^\vep\otimes\P'(d\t_k,d\omega'_k)F\left(\theta^{(\vep)}(\t_1,\omega'_1),\ldots,\theta^{(\vep)}(\t_k,\omega'_k)\right)\label{arbreGW2}.
\end{eqnarray}
Note that for every permutation $\vph$ of the set $\{1,\ldots,k\}$, $(\tau_{\vph(1)}\theta^{(\vep)},\ldots,\tau_{\vph(k)}\theta^{(\vep)})$ and $(\tau_1\theta^{(\vep)},\ldots,\tau_k\theta^{(\vep)})$ have the same distribution under $\tt^\vep\otimes\P'$. Then, (\ref{arbreGW2}) means that the law of $\theta^{(\vep)}$ under $\tt^\vep\otimes\P'$ satisfies the branching property of the Galton-Watson trees. This completes the proof of the desired result.\cq

Recall that $\sum_{i\in I}\dl_{\t_i}$ is a Poisson point measure on $\T$ with intensity $\tt$. Let us now set, for every $t,h>0$,
$$\z(t,t+h)=\sum_{i\in I}Z(t,t+h)(\t_i).$$
For every $\vep>0$, we define a process $\x^\vep=(\x^\vep_k,k\geq0)$ on $(\Omega,\P)$ by the formula
$$\x^\vep_k=\z(k\vep,(k+1)\vep),\;k\geq0.$$

\begin{proposition}\label{GaltonWatson}
For every $\vep>0$, the process $\x^\vep$ is a Galton-Watson process whose initial distribution is the Poisson distribution with parameter $v(\vep)$ and whose offspring distribution is $\mu_\vep$.
\end{proposition}

\proof We first observe that $\x^\vep_0=\n(\h(\t)>\vep)$ is Poisson with parameter \linebreak $\tt(\h(\t)>\vep)=v(\vep)$. Then let $p$ be a positive integer. We know from a classical property of Poisson measures that, under the probability measure $\P$ and conditionally on the event $\{\x^\vep_0=p\}$, the atoms of $\n$ that belong to the set $\T^\vep$ are distributed as $p$ i.i.d. variables with distribution $\tt^\vep$. Furthermore, it follows from Lemma \ref{arbreGW} that under $\tt^\vep$, the process $(Z(k\vep,(k+1)\vep))_{k\geq0}$ is a Galton-Watson process started at one with offspring distribution $\mu_\vep$. This completes the proof.\cq

As a consequence, we get the next proposition, which we will use throughout this work.

\begin{proposition}\label{souscritique}
For every $t>0$ and $h>0$, we have $\tt(Z(t,t+h))\leq v(h)$.
\end{proposition}

\proof Since compact real trees have finite height, the Galton-Watson process $\x^\vep$ dies out $\P$ a.s. This implies that $\mu_\vep$ is critical or subcritical so that $(\x^\vep_k,k\geq0)$ is a supermartingale. Let $t,h>0$. We can find $\vep>0$ and $k\in\N$ such that $t=k\vep$ and $\vep\leq h$. Thus we have, 
\be\label{inegsouscritique}
\tt(Z(t,t+\vep))=\tt(Z(k\vep,(k+1)\vep))=\E(\x^\vep_k)\leq\E(\x^\vep_0)=v(\vep).
\ee  
Using Lemma \ref{binomiale} and (\ref{inegsouscritique}), we get
$$\tt(Z(t,t+h))=\tt\left(Z(t,t+\vep)\frac{v(h)}{v(\vep)}\right)\leq v(h).$$\cq

\subsubsection{A local time process}

\begin{proposition}\label{defL_t}
For every $t\geq0$, there exists a random variable $L_t$ on the space $\T$ such that $\tt$ a.e.,
$$\frac{Z(t,t+h)}{v(h)}\build{\longrightarrow}_{h\to 0}^{} L_t.$$
\end{proposition}

\proof Let us start with the case $t=0$. As $Z(0,h)=\ind{\{\h(\t)>h\}}$ for every $h>0$, Lemma \ref{vcontinue} gives $v(h)^{-1}Z(0,h)\la0$ $\tt$ a.e. as $h\la0$, so we set $L_0=0$.

Let us now fix $t>0$. Thanks to Lemma \ref{vcontinue}, we can define a decreasing sequence $(\vep_n)_{n\geq1}$ by the condition $v(\vep_n)=n^4$ for every $n\geq1$. We claim that there exists a random variable $L_t$ on the space $\T$ such that, $\tt$ a.e.,
\begin{equation}\label{convL_tps}
\frac{Z(t,t+\vep_n)}{n^4}\build{\longrightarrow}_{n\to\infty}^{}L_t.
\end{equation}
Indeed, using Lemma \ref{binomiale}, we have, for every $n\geq1$, 
\begin{eqnarray}
&&\hspace{-1cm}\tt^t\left(\left|\frac{Z(t,t+\vep_n)}{n^4}-\frac{Z(t,t+\vep_{n+1})}{(n+1)^4}\right|^2\right)\nonumber\\
&=&\tt^t\left(\frac{1}{n^{8}}\tt^t\left(\left|Z(t,t+\vep_n)-\frac{n^4}{(n+1)^4}Z(t,t+\vep_{n+1})\right|^2 \Big|\; Z(t,t+\vep_{n+1})\right)\right)\nonumber\\
&\leq&\tt^t\left(\frac{Z(t,t+\vep_{n+1})}{4n^{8}}\right)\nonumber\\
&\leq&\frac{(n+1)^4}{4v(t)n^{8}},\label{Palm1}
\end{eqnarray} 
where the last bound follows from Proposition \ref{souscritique} and the definition of $\vep_{n+1}$. Thanks to the Cauchy-Schwarz inequality, we get
\begin{equation}\label{normeL1}
\tt^t\left(\left| \frac{Z(t,t+\vep_n)}{n^4}-\frac{Z(t,t+\vep_{n+1})}{(n+1)^4}\right|\right)\leq \frac{(n+1)^2}{2n^4\sqrt{v(t)}}\leq\frac{2}{n^2\sqrt{v(t)}}.
\end{equation}
The bound (\ref{normeL1}) implies 
$$\tt\left(\sum_{n=1}^{\infty}\left| \frac{Z(t,t+\vep_n)}{n^4}-\frac{Z(t,t+\vep_{n+1})}{(n+1)^4}\right|\right)<\infty.$$
In particular, $\tt$ a.e.,
$$\sum_{n=1}^{\infty}\left| \frac{Z(t,t+\vep_n)}{n^4}-\frac{Z(t,t+\vep_{n+1})}{(n+1)^4}\right|<\infty.$$
Our claim (\ref{convL_tps}) follows.

For every $h\in(0,\vep_1]$, we can find $n\geq1$ such that $\vep_{n+1}\leq h\leq\vep_n$. Then, we have $Z(t,t+\vep_n)\leq Z(t,t+h)\leq Z(t,t+\vep_{n+1})$ $\tt$ a.e., and $n^4\leq v(h)\leq (n+1)^4$ so that
$$\frac{Z(t,t+\vep_n)}{(n+1)^4}\leq\frac{Z(t,t+h)}{v(h)}\leq\frac{Z(t,t+\vep_{n+1})}{n^4}.$$
We then deduce from (\ref{convL_tps}) that $\tt$ a.e.,
$$\frac{Z(t,t+h)}{v(h)}\build{\longrightarrow}_{h\to0}^{}L_t$$
which completes the proof.\cq

\begin{definition}
We define a process $\l=(\l_t,t\geq0)$ on $(\Omega,\P)$ by setting $\l_0=1$ and for every $t>0$,
$$\l_t=\sum_{i\in I}L_t(\t_i).$$
\end{definition}
Notice that $L_t(\t)=0$ if $\h(\t)\leq t$ so that the above sum is finite a.s.
 
\begin{corollary}\label{theoreml_t}
For every $t\geq0$, we have $\P$ a.s.
$$\frac{\z(t,t+h)}{v(h)}\build{\longrightarrow}_{h\to0}^{}\l_t.$$
Moreover, this convergence holds in $\L^1(\P)$ uniformly in $t\in[0,\infty)$.
\end{corollary}

\proof The first assertion is an immediate consequence of Proposition \ref{defL_t}. Let us focus on the second assertion. From Lemma \ref{binomiale}, $\tt(n^{-4}Z(t,t+\vep_n)-(n+1)^{-4}Z(t,t+\vep_{n+1}))=0$ for every $t\geq0$ and $n\geq1$. Thus, from the second moment formula for Poisson measures, we get, for every $t\geq0$ and $n\geq1$,
$$\E\left(\left(\frac{\z(t,t+\vep_n)}{n^4}-\frac{\z(t,t+\vep_{n+1})}{(n+1)^4}\right)^2\right)=\tt\left(\left(\frac{Z(t,t+\vep_n)}{n^4}-\frac{Z(t,t+\vep_{n+1})}{(n+1)^4}\right)^2\right).$$
Now, we have
$$\tt\left(\left(\frac{Z(0,\vep_n)}{n^4}-\frac{Z(0,\vep_{n+1})}{(n+1)^4}\right)^2\right)=\tt\left(\left(\frac{\ind{\{\h(\t)>\vep_n\}}}{n^4}-\frac{\ind{\{\h(\t)>\vep_{n+1}\}}}{(n+1)^4}\right)^2\right)=\frac{1}{n^4}-\frac{1}{(n+1)^4}$$
and for every $t>0$, thanks to the bound (\ref{Palm1}),
$$\tt\left(\left(\frac{Z(t,t+\vep_n)}{n^4}-\frac{Z(t,t+\vep_{n+1})}{(n+1)^4}\right)^2\right)\leq\frac{(n+1)^4}{4n^8}.$$
So for every $t\geq0$ and $n\geq1$, we have from the Cauchy-Schwarz inequality
\be\label{Palm2}
\E\left(\left|\frac{\z(t,t+\vep_n)}{n^4}-\frac{\z(t,t+\vep_{n+1})}{(n+1)^4}\right|\right)\leq\frac{(n+1)^2}{n^4}.
\ee
Then $n^{-4}\z(t,t+\vep_n)\la\l_t$ in $\L^1$ as $n\la\infty$ and, for every $n\geq2$,
$$\E\left(\left|\frac{\z(t,t+\vep_n)}{n^4}-\l_t\right|\right)\leq\sum_{k=n}^\infty\frac{(k+1)^2}{k^4}\leq\sum_{k=n}^\infty\frac{4}{k^2}\leq\frac{8}{n}.$$
In the same way as in the proof of (\ref{Palm2}), we have the following inequality: If $h\in(0,\vep_1]$, $t\geq0$ and $n$ is a positive integer such that $\vep_{n+1}\leq h\leq\vep_n$,
$$\E\left(\left|\frac{\z(t,t+\vep_n)}{n^4}-\frac{\z(t,t+h)}{v(h)}\right|\right)\leq\frac{\sqrt{v(h)}}{n^4}\leq\frac{16}{\sqrt{v(h)}}.$$ 
Then, for every $h\in(0,\vep_2]$ and $t\geq0$, we get 
$$\E\left(\left|\frac{\z(t,t+h)}{v(h)}-\l_t\right|\right)\leq16\left(v(h)^{-1/2}+v(h)^{-1/4}\right),$$
which completes the proof.\cq

We will now establish a regularity property of the process $(\l_t,t\geq0)$. 

\begin{proposition}\label{modif}
The process $(\l_t,t\geq0)$ admits a modification, denoted by $(\wt{\l}_t,t\geq0)$, which is right-continuous with left-limits, and which has no fixed discontinuities.
\end{proposition} 

\proof We start with two lemmas.

\begin{lemma}\label{E(l_t)continue}
There exists $\lm\geq0$ such that $\E(\l_t)=e^{-\lm t}$ for every $t\geq0$.
\end{lemma}

\proof We claim that the function $t\in[0,+\infty)\lmt\E(\l_t)$ is multiplicative, meaning that for every $t,s\geq0$, $\E(\l_{t+s})=\E(\l_t)\E(\l_s)$. As $\l_0=1$ by definition, $\E(\l_0)=1$. Let $t,s>0$ and $0<h<s$. Let us denote by $\t^1,\ldots,\t^{Z(t,t+h)}$ the subtrees of $\t$ above level $t$ with height greater than $h$. Then, using the regenerative property, we can write
$$\tt(Z(t+s,t+s+h))=\tt\left(\sum_{i=1}^{Z(t,t+h)}Z(s,s+h)(\t^i)\right)=\tt\left(Z(t,t+h))\tt^h(Z(s,s+h)\right),$$
which implies
\be\label{mult}
\E(\z(t+s,t+s+h))=\E(\z(t,t+h))\E\left(\frac{\z(s,s+h)}{v(h)}\right).
\ee
 Thus, dividing by $v(h)$ and letting $h\la0$ in (\ref{mult}), we get our claim from Corollary \ref{theoreml_t}. Moreover, thanks to Proposition \ref{souscritique} and Corollary \ref{theoreml_t}, we know that $\E(\l_t)\leq1$ for every $t\geq0$. Then, we obtain in particular that the function $t\in[0,\infty)\lmt\E(\l_t)$ is nonincreasing.  

To complete the proof, we have to check that $\E(\l_t)>0$ for every  $t>0$. If we assume that $\E(\l_t)=0$ for some $t>0$ then $L_t=0$, $\tt$ a.e. Let $s,h>0$ such that $0<h<s$. With the same notation as in the beginning of the proof, we can write 
\begin{eqnarray}
\tt(\h(\t)>t+s)&=&\tt\left(\exists i\in\{1,\ldots,Z(t,t+h)\}:\;\h(\t^i)>s\right)\nonumber\\
&=&\tt\left(1-\left(1-\frac{v(s)}{v(h)}\right)^{Z(t,t+h)}\right).\label{E(l_t)>0}
\end{eqnarray}  
Now, thanks to Proposition \ref{defL_t}, $\tt$ a.e.,
$$\left(1-\frac{v(s)}{v(h)}\right)^{Z(t,t+h)}\build{\longrightarrow}_{h\to0}^{}\;\exp(-L_tv(s))=1.$$ 
Moreover, $\tt$ a.e.,
$$1-\left(1-\frac{v(s)}{v(h)}\right)^{Z(t,t+h)}\leq\ind{\{\h(\t)>t\}}.$$
Then, using dominated convergence in (\ref{E(l_t)>0}) as $h\la0$, we obtain $\tt(\h(\t)>t+s)=0$ which contradicts the assumptions of Theorem \ref{casinfini}.\cq

\begin{lemma}\label{surmart}
Let us denote by $D=\lbrace k2^{-n},k\geq1,n\geq0\rbrace$ the set of positive dyadic numbers and define $\g_t=\sg(\l_s,s\in D,s\leq t)$ for every $t\in D$. Then $(\l_t,t\in D)$ is a nonnegative supermartingale with respect to the filtration $(\g_t,t\in D)$. 
\end{lemma}

\proof Let $p$ be a positive integer, let $s_1,\ldots,s_p,s,t\in D$ such that $s_1<\ldots<s_p\leq s<t$ and let $f:\R^p\longrightarrow\R_+$ be a bounded continuous function. We can find a positive integer $n$ such that $2^nt$, $2^ns$, and $2^ns_i$ for $i\in\lbrace 1,\ldots,p\rbrace$ are nonnegative integers. The process $\x^{2^{-n}}$ is a subcritical Galton-Watson process, so  
$$\E(\x^{2^{-n}}_{2^nt}f(\x^{2^{-n}}_{2^ns_1},\ldots,\x^{2^{-n}}_{2^ns_p}))\leq\E(\x^{2^{-n}}_{2^ns}f(\x^{2^{-n}}_{2^ns_1},\ldots,\x^{2^{-n}}_{2^ns_p})).$$ 
Therefore we have also,
\begin{eqnarray}
&&\hspace{-1cm}\E\left(\frac{\z(t,t+2^{-n})}{v(2^{-n})}f\left(\frac{\z(s_1,s_1+2^{-n})}{v(2^{-n})},\ldots,\frac{\z(s_p,s_p+2^{-n})}{v(2^{-n})}\right)\right)\nonumber\\
&\leq&\E\left(\frac{\z(s,s+2^{-n})}{v(2^{-n})}f\left(\frac{\z(s_1,s_1+2^{-n})}{v(2^{-n})},\ldots,\frac{\z(s_p,s_p+2^{-n})}{v(2^{-n})}\right)\right)\label{surmart2}.
\end{eqnarray}
We can then use Corollary \ref{theoreml_t} to obtain $\E(\l_tf(\l_{s_1},\ldots,\l_{s_p}))\leq\E(\l_sf(\l_{s_1},\ldots,\l_{s_p}))$. \cq

We now complete the proof of Proposition \ref{modif}. Let us set, for every $t\geq0$,
$$\wt{\g}_t=\bigcap_{s>t,s\in D}\g_s.$$  
From Lemma \ref{surmart} and classical results on supermartingales, we can define a right-\linebreak continuous supermartingale $(\wt{\l}_t,t\geq0)$ with respect to the filtration $(\wt{\g}_t,t\geq0)$ by setting, for every $t\geq0$,
\be\label{defwtl_t}
\wt{\l}_t=\build{\lim}_{s\downarrow t,s\in D}^{}\l_s,
\ee
where the limit holds $\P$ a.s. and in $\L^1$ (see e.g. Chapter VI in \cite{DM} for more details). We claim that $(\wt{\l}_t,t\geq0)$ is a c\`adl\`ag modification of $(\l_t,t\geq0)$ with no fixed discontinuities.

We first prove that $(\wt{\l}_t,t\geq0)$ is a modification of $(\l_t,t\geq0)$. For every $t\geq0$ and every sequence $(s_n)_{n\geq0}$ in $D$  such that $s_n\da t$ as $n\ua\infty$, we have thanks to (\ref{defwtl_t}) and Lemma \ref{E(l_t)continue},
$$\E(\wt{\l}_t)=\lim_{n\to\infty}\E(\l_{s_n})=\E(\l_t).$$
Let us now show that for every $t\geq0$, $\l_t\leq\wt{\l_t}$ $\P$ a.s. Let $\al,\vep>0$ and  $\dl\in(0,1)$. Thanks to Corollary \ref{theoreml_t}, we can find $h_0>0$ such that for every $h\in(0,h_0)$ and $n\geq0$, 
$$\E\left(\left|\frac{\z(t,t+h)}{v(h)}-\l_t\right|\right)\leq\vep\al\;\;{\rm and}\;\;\E\left(\left|\frac{\z(s_n,s_n+h)}{v(h)}-\l_{s_n}\right|\right)\leq\vep\al.$$
We choose $h\in(0,h_0)$ and $n_0\geq0$ such that $s_n-t+h\leq h_0$ and $v(h)\leq(1+\dl)v(s_n-t+h)$ for every $n\geq n_0$. We notice that $\z(t,s_n+h)\leq\z(s_n,s_n+h)$ so that, for every $n\geq n_0$,
\begin{eqnarray*}
&&\hspace{-1cm}\P(\l_t>(1+\dl)\l_{s_n}+\vep)\\
&\leq&\P\left(\l_t-\frac{\z(t,s_n+h)}{v(s_n-t+h)}>(1+\dl)\l_{s_n}-(1+\dl)\frac{\z(s_n,s_n+h)}{v(h)}+\vep\right)\\
&\leq&2\vep^{-1}\E\left(\left|\frac{\z(t,s_n+h)}{v(s_n-t+h)}-\l_t\right|\right)+2\vep^{-1}(1+\dl)\E\left(\left|\frac{\z(s_n,s_n+h)}{v(h)}-\l_{s_n}\right|\right)\\
&\leq&6\alpha.
\end{eqnarray*}
We have thus shown that
\begin{equation}\label{limmart}
\P(\l_t>(1+\dl)\l_{s_n}+\vep)\build{\longrightarrow}_{n\to\infty}^{}0.
\end{equation}
So, $\P(\l_t-(1+\dl)\wt{\l}_t>\vep)=0$ for every $\vep>0$, implying that $\l_t\leq (1+\dl)\wt{\l}_t$, $\P$ a.s. This leads us to the claim $\l_t\leq\wt{\l}_t$ a.s. Since we saw that $\E(\l_t)=\E(\wt{\l}_t)$, we have $\l_t=\wt{\l}_t$ $\P$ a.s. for every $t\geq0$.

Now, $(\wt{\l}_t,t\geq0)$ is a right-continuous supermartingale. Thus, $(\wt{\l}_t,t\geq0)$ is also left-limited and we have $\E(\wt{\l}_t)\leq\E(\wt{\l}_{t-})$ for every $t>0$. Moreover, we can prove in the same way as we did for (\ref{limmart}) that, for every $t>0$ and every sequence $(s_n,n\geq0)$ in $D$ such that $s_n\ua t$ as $n\ua\infty$,
$$\P(\wt{\l}_{s_n}>(1+\dl)\wt{\l}_t+\vep)\build{\longrightarrow}_{n\to\infty}^{}0,$$
implying that $\wt{\l}_{t-}\leq\wt{\l}_t$, $\P$ a.s. So, $\l_t=\l_{t-}$ $\P$ a.s. for every $t>0$ meaning that $(\wt{\l}_t,t\geq0)$ has no fixed discontinuities. \cq

From now on, to simplify notation, we replace $(\l_t,t\geq0)$ by its c\`adl\`ag modification $(\wt{\l}_t,t\geq0)$.

\subsubsection{The CSBP}

We will prove that the suitably rescaled family of Galton-Watson processes $(\x^\vep)_{\vep>0}$ converges to the local time $\l$.

Thanks to Lemma \ref{vcontinue}, we can define a sequence $(\eta_n)_{n\geq1}$ by the condition $v(\eta_n)=n$ for every $n\geq1$. We set $m_n=[\eta_n^{-1}]$ where $[x]$ denotes the integer part of $x$. We recall from Proposition \ref{GaltonWatson} that $\x^{\eta_n}$ is a Galton-Watson process on $(\Omega,\P)$ whose initial distribution is the Poisson distribution with parameter $n$. For every $n\geq1$, we define a process $\y^n=(\y^n_t,t\geq0)$ on $(\Omega,\P)$ by the following formula,
$$\y^n_t=n^{-1}\x^{\eta_n}_{[m_nt]},\;t\geq0.$$
 
\begin{proposition}\label{convl_t1}
For every $t\geq0$, $\y^n_t\la\l_t$ in probability as $n\la\infty$.
\end{proposition} 

\proof The result for $t=0$ is a consequence of the definition of $\l_0$ together with simple estimates for the Poisson distribution. Let $t,\dl>0$. We can write
\begin{eqnarray*}
\P(\vert\y^n_t-\l_t\vert>2\dl)
&\leq&\P\left(\left|\y^n_t-\l_{\eta_n[m_nt]}\right|>\dl\right)+\P\left(\left|\l_{\eta_n[m_nt]}-\l_t\right|>\dl\right)\\
&\leq&\dl^{-1}\E\left(\left|\y^n_t-\l_{\eta_n[m_nt]}\right|\right)+\P\left(\left|\l_{\eta_n[m_nt]}-\l_t\right|>\dl\right).
\end{eqnarray*}
Now, Corollary \ref{theoreml_t} and Proposition \ref{modif} imply respectively that
$$\E\left(\left|\y^n_t-\l_{\eta_n[m_nt]}\right|\right)\build{\longrightarrow}_{n\to\infty}^{}0\;\;{\rm and}\;\;\P\left(\left|\l_{\eta_n[m_nt]}-\l_t\right|>\dl\right)\build{\longrightarrow}_{n\to\infty}^{}0,$$
which completes the proof.\cq

\begin{corollary}
For every $t\geq0$, the law of $\y^n_t$ under $\P(\cdot\mid\x^{\eta_n}_0=n)$ converges weakly to the law of $\l_t$ under $\P$ as $n\la\infty$. 
\end{corollary}

\proof For positive integers $n$ and $k$, we denote by $\mu^n$ the offspring distribution of the Galton-Watson process $\x^{\eta_n}$, by $f^n$ the generating function of $\mu^n$ and by $f_k^n$ the k-th iterate $f_k^n=f^n\circ\ldots\circ f^n$ of $f^n$. Let $\lm>0$ and $t\geq0$. We  have,
$$\E(\exp(-\lm \y^n_t))=\sum_{p=0}^{\infty}e^{-n}\frac{n^p}{p!}\left(f^n_{[m_nt]}\left(e^{-\lm/n}\right)\right)^p=\exp\left(-n\left(1-f^n_{[m_nt]}\left(e^{-\lm/n}\right)\right)\right).$$
From Proposition \ref{convl_t1}, it holds that
$$\exp\left(-n\left(1-f^n_{[m_nt]}\left(e^{-\lm/n}\right)\right)\right)\build{\longrightarrow}_{n\to\infty}^{}\E(\exp(-\lm \l_t)).$$
Let us set $u(t,\lm)=-\log(\E[\exp(-\lm\l_t)])$. It follows that,
$$n\left(1-f^n_{[m_nt]}\left(e^{-\lm/n}\right)\right)\build{\la}_{n\to\infty}^{}u(t,\lm).$$
Thus, we obtain,
$$\E\left(\exp(-\lm \y^n_t)\mid\x^{\eta_n}_0=n\right)=\left(f^n_{[m_nt]}\left(e^{-\lm/n}\right)\right)^n\build{\longrightarrow}_{n\to\infty}^{}\exp(-u(t,\lm))=\E[\exp(-\lm \l_t)].$$\cq

At this point, we can use Theorem \ref{Grimvall} to assert that $(\l_t,t\geq0)$ is a CSBP and that the law of $(\y^n_t,t\geq0)$ under the probability measure $\P(\cdot\mid\x^{\eta_n}_0=n)$ converges to the law of $(\l_t,t\geq0)$ as $n\la\infty$ in the space of probability measures on the Skorokhod space $\D(\R_+)$. To verify the assumptions of Theorem \ref{Grimvall}, we need to check that there exists $\dl>0$ such that $\P(\l_\dl>0)>0$. This is obvious from Lemma \ref{E(l_t)continue}.

\subsection{Identification of the measure $\tt$}

In the previous section, we have constructed from $\tt$ a CSBP $\l$, which becomes extinct almost surely. We denote by $\psi$ the associated  branching mechanism. We can consider the $\sg$-finite measure $\tt_\psi$, which is the law of the L\'evy tree associated with $\l$. Our goal is to show that the measures $\tt$ and $\tt_\psi$ coincide.

Recall that $\mu^n$ denotes the offspring distribution of the Galton-Watson process $\x^{\eta_n}$.
 
\begin{lemma}\label{DuLG2}
For every $a>0$, the law of the $\R$-tree $\eta_n\t^\theta$ under $\Pi_{\mu^n}(\cdot\mid\h(\theta)\geq [am_n])$ converges as $n\la\infty$ to the probability measure $\tt_\psi(\cdot\mid\h(\t)>a)$ in the sense of weak convergence of measures in the space $\T$. 
\end{lemma}

\proof We first check that, for every $\dl>0$,
\be\label{limDuLG}
\liminf_{n\rightarrow\infty} \P(\y^n_\dl=0)>0.
\ee
Indeed, we have
$$\P(\y^n_\dl=0)=\P(\n(\h(\t)>([m_n\dl]+1)\eta_n)=0)=\exp(-v(([m_n\dl]+1)\eta_n)).$$
As $v$ is continuous, it follows that $\P(\y^n_\dl=0)\la\exp(-v(\dl))$ as $n\la\infty$ implying (\ref{limDuLG}).

We recall that the law of $\y^n$ under the probability measure $\P(\cdot\mid\x^{\eta_n}_0=n)$ converges to the law of $(\l_t,t\geq0)$. Then, thanks to (\ref{limDuLG}), we can apply Theorem \ref{DuLG} to get that, for every $a>0$, the law of the $\R$-tree $m_n^{-1}\t^\theta$ under $\Pi_{\mu^n}(\cdot\mid\h(\theta)\geq [am_n])$ converges to the probability measure $\tt_\psi(\cdot\mid\h(\t)>a)$ in the sense of weak convergence of measures in the space $\T$. As $m_n^{-1}\eta_n\la1$ as $n\la\infty$, we get the desired result.\cq

We can now complete the proof of Theorem \ref{casinfini}. Indeed, thanks to Lemmas \ref{discretisation} and \ref{arbreGW}, we can construct on the same probability space $(\OO,\PP)$, a sequence of $\T$-valued random variables $(\TT_n)_{n\geq1}$ distributed according to $\tt(\cdot\mid\h(\t)>([am_n]+1)\eta_n)$ and a sequence of $\a$-valued random variables $(\ttt_n)_{n\geq1}$ distributed according to $\Pi_{\mu^n}(\cdot\mid\h(\theta)\geq[am_n])$ such that for every $n\geq1$, $\PP$ a.s.,
$$d_{GH}(\TT_n,\eta_n\t^{\ttt_n})\leq4\eta_n.$$
Then, using Lemma \ref{DuLG2}, we have $\tt(\cdot\mid\h(\t)>([am_n]+1)\eta_n)\la\tt_\psi(\cdot\mid\h(\t)>a)$ as $n\la\infty$ in the sense of weak convergence of measures on the space $\T$. So we get  
$$\tt(\cdot\;\vert\h(\t)>a)=\tt_\psi(\cdot\;\vert\h(\t)>a)$$ 
for every $a>0$, and thus $\tt=\tt_\psi$.

\section{Proof of Theorem \ref{casfini}}\label{sec4}

Let $\tt$ be a probability measure on $(\T,d_{GH})$ satisfying the assumptions of Theorem \ref{casfini}.

In this case, we define $v:[0,\infty)\la(0,\infty)$ by $v(t)=\tt(\h(\t)>t)$ for every $t\geq0$. Note that $v(0)=1$ is well defined here. For every $t>0$, we denote by $\tt^t$ the probability measure $\tt(\cdot\mid\h(\t)>t)$. The following two results are proved in a similar way to Lemmas \ref{vcontinue} and \ref{binomiale}.

\begin{lemma}\label{vcontinue2}
The function $v$ is nonincreasing, continuous and goes to $0$ as $t\la\infty$.
\end{lemma}

\begin{lemma}\label{binomiale2}
For every $t>0$ and $0<a<b$, the conditional law of the random variable $Z(t,t+b)$, under the probability measure $\tt^t$ and given $Z(t,t+a)$, is a binomial distribution with parameters $Z(t,t+a)$ and $v(b)/v(a)$. 
\end{lemma}

\subsection{The DSBP derived from $\tt$}

We will follow the same strategy as in section \ref{sec3} but instead of a CSBP we will now construct an integer-valued branching process.

\subsubsection{A family of Galton-Watson trees}

We recall that $\mu_\vep$ denotes the law of $Z(\vep,2\vep)$ under the probability measure $\tt^\vep$, and that $(\theta_\xi,\xi\in\A)$ is a sequence of independent $\a$-valued random variables defined on a probability space $(\Omega',\P')$ such that for every $\xi\in\A$, $\theta_\xi$ is distributed uniformly over $\proj^{-1}(\xi)$. The following lemma is proved in the same way as Lemma \ref{arbreGW}. 

\begin{lemma}\label{arbreGWbis}
Let us define for every $\vep>0$, a mapping $\theta^{(\vep)}$ from $\T^{(\vep)}\times\Omega'$ into $\a$ by 
$$\theta^{(\vep)}(\t,\omega)=\theta_{\xi^\vep(\t)}(\omega).$$
Then for every positive integer $p$, the law of the random variable $\theta^{(\vep)}$ under the probability measure $\tt^{p\vep}\otimes\P'$ is $\Pi_{\mu_\vep}(\cdot\mid\h(\theta)\geq p-1)$.
\end{lemma} 

For every $\vep>0$, we define a process $X^\vep=(X^\vep_k,k\geq0)$ on $\T$ by the formula
$$X^\vep_k=Z(k\vep,(k+1)\vep),\;k\geq0.$$ 
We show in the same way as Proposition \ref{GaltonWatson} and Proposition \ref{souscritique} the following two results.

\begin{proposition}\label{GaltonWatson2}
For every $\vep>0$, the process $X^\vep$ is under $\tt$ a Galton-Watson process whose initial distribution is the Bernoulli distribution with parameter $v(\vep)$ and whose offspring distribution is $\mu_\vep$.
\end{proposition}

\begin{proposition}\label{souscritique2}
For every $t>0$ and $h>0$, we have $\tt(Z(t,t+h))\leq v(h)\leq1$.
\end{proposition}

The next proposition however is particular to the finite case and will be useful in the rest of this section.

\begin{proposition}\label{dl_1}
The family of probability measures $(\mu_\vep)_{\vep>0}$ converges to the Dirac measure $\dl_1$ as $\vep\to0$. In other words,
$$\tt^\vep(Z(\vep,2\vep)=1)\build{\la}_{\vep\to0}^{}1.$$
\end{proposition}

\proof We first note that 
$$2\tt^\vep(Z(\vep,2\vep)\geq1)-\tt^\vep(Z(\vep,2\vep))\leq\tt^\vep(Z(\vep,2\vep)=1)\leq\tt^\vep(Z(\vep,2\vep)\geq1).$$
Moreover, $\tt^\vep(Z(\vep,2\vep)\geq 1)=\tt^\vep(\h(\t)>2\vep)=v(2\vep)/v(\vep)$ and $\tt^\vep(Z(\vep,2\vep))\leq1$. So,  
\be\label{inegdl_1}
\frac{2v(2\vep)}{v(\vep)}-1\leq\tt^\vep(Z(\vep,2\vep)=1)\leq\frac{v(2\vep)}{v(\vep)}.
\ee
We let $\vep\la0$ in (\ref{inegdl_1}) and we use Lemma \ref{vcontinue2} to obtain the desired result.\cq

\subsubsection{Construction of the DSBP}

\begin{proposition}\label{defL_t2}
For every $t\geq0$, there exists an integer-valued random variable $L_t$ on the space $\T$ such that $\tt(L_t)\leq1$ and $\tt$ a.s.,
$$Z(t,t+h)\build{\ua}_{h\da0}^{}L_t.$$
\end{proposition} 

\proof Let $t\geq0$. The function $h\in(0,\infty)\longmapsto Z(t,t+h)\in\Z_+$ is nonincreasing so that there exists a random variable $L_t$ with values in $\Z_+\cup\{\infty\}$ such that, $\tt$ a.s., 
$$Z(t,t+h)\build{\ua}_{h\da0}^{}L_t.$$ 
Thanks to the monotone convergence theorem, we have
$$\tt(Z(t,t+h))\build{\la}_{h\to0}^{}\tt(L_t).$$ 
Now, by Proposition \ref{souscritique2}, $\tt(Z(t,t+h))\leq1$ for every $h>0$. Then, $\tt(L_t)\leq1$ which implies in particular that $L_t<\infty$ $\tt$ a.s.\cq

\begin{proposition}\label{conv2L_t} For every $t>0$, the following two convergences hold $\tt$ a.s., 
\begin{eqnarray}
Z(t-h,t)&\build{\ua}_{h\da0}^{}&L_t\label{conv2L_t1},\\
Z(t-h,t+h)&\build{\ua}_{h\da0}^{}&L_t\label{conv2L_t2}.
\end{eqnarray}
\end{proposition}

\proof Let $t>0$ be fixed throughout this proof. By the same arguments as in the proof of Proposition \ref{defL_t2}, we can find a $\Z_+$-valued random variable $\overline{L}_t$ such that $\tt(\overline{L}_t)\leq1$ and $Z(t-h,t)\ua\overline{L}_t$ as $h\da0$, $\tt$ a.s. If $h\in(0,t)$, we write $\t^1,\ldots,\t^{Z(t-h,t)}$ for the subtrees of $\t$ above level $t-h$ with height greater than $h$. Then, from the regenerative property, 
\begin{eqnarray}
&&\hspace{-1cm}\tt\left(\vert Z(t,t+h)-Z(t-h,t)\vert\geq1\right)\nonumber\\
&=& \tt\Big(\tt\Big(\Big|\sum_{i=1}^{Z(t-h,t)}(Z(h,2h)(\t^{i})-1)\Big|\geq1\;\Big|\; Z(t-h,t)\Big)\Big)\nonumber\\
&\leq&\tt\left(\tt\left(|Z(h,2h)(\t^i)-1|\geq1\;\;{\rm for\; some}\;\;i\in\{1,\ldots,Z(t-h,t)\}\mid Z(t-h,t)\right)\right)\nonumber\\
&\leq&\tt\left(Z(t-h,t)\tt^h\left(\vert Z(h,2h)-1\vert\geq1\right)\right)\label{eqconv2L_t1}.
\end{eqnarray}
Since $Z(t-h,t)\tt^h\left(\vert Z(h,2h)-1\vert\geq1\right)\leq\overline{L}_t$ $\tt$ a.s., Proposition \ref{dl_1} and the dominated convergence theorem imply that the right-hand side of (\ref{eqconv2L_t1}) goes to $0$ as $h\la0$. Thus $L_t=\overline{L}_t$ $\tt$ a.s.

Likewise, there exists a random variable $\wh{L}_t$ with values in $\Z_+$ such that, $\tt$ a.s., \linebreak $Z(t-h,t+h)\ua\wh{L}_t$ as $h\da0$. Let us now notice that, for every $h>0$, $\tt$ a.s., \linebreak $Z(t-h,t+h)\leq Z(t-h,t)$.  Moreover, thanks to Lemma \ref{binomiale2}, we have
\be\label{eqconv2L_t2}
\tt(Z(t-h,t)\geq Z(t-h,t+h)+1)=1-\tt\left(\left(\frac{v(2h)}{v(h)}\right)^{Z(t-h,t)}\right)\geq1-\tt\left(\left(\frac{v(2h)}{v(h)}\right)^{L_t}\right).
\ee
The right-hand side of (\ref{eqconv2L_t2}) tends to $0$ as $h\la0$. So $L_t=\wh{L}_t$ $\tt$ a.s. \cq

We will now establish a regularity property of the process $(L_t,t\geq0)$. 

\begin{proposition}\label{L_tcadlag}
The process $(L_t,t\geq0)$ admits a modification which is right-continuous with left limits, and which has no fixed discontinuities.
\end{proposition}

\proof We start the proof with three lemmas. The first one in proved in a similar but easier way as Lemma \ref{E(l_t)continue}.

\begin{lemma}\label{tt(L_t)continue}
There exists $\lm\geq0$ such that $\tt(L_t)=e^{-\lm t}$ for every $t\geq0$.
\end{lemma}

For every $n\geq1$ and every $t\geq0$ we set $Y^n_t=X^{1/n}_{[nt]}$. 

\begin{lemma}\label{convDSBP}
For every $t\geq0$, $Y^n_t\la L_t$ as $n\la\infty$, $\tt$ a.s.
\end{lemma}
This lemma is an immediate consequence of Proposition \ref{conv2L_t}.

\begin{lemma}\label{Lsurmart}
Let us define $G_t=\sg(L_s,s\leq t)$ for every $t\geq0$. Then $(L_t,t\geq0)$ is a nonnegative supermartingale with respect to the filtration $(G_t,t\geq0)$.
\end{lemma}

\proof Let $s,t,s_1,\ldots,s_p\geq0$ such that $0\leq s_1\leq\ldots\leq s_p\leq s<t$ and let $f:\R^p\rightarrow\R_+$ be a bounded measurable function. For every $n\geq1$, the offspring distribution $\mu_{1/n}$ is critical or subcritical so that $(X^{1/n}_k,k\geq0)$ is a supermartingale. Thus we have 
$$\tt\left(X^{1/n}_{[nt]}f\left(X^{1/n}_{[ns_1]},\ldots,X^{1/n}_{[ns_p]}\right)\right)\leq\tt\left(X^{1/n}_{[ns]}f\left(X^{1/n}_{[ns_1]},\ldots,X^{1/n}_{[ns_p]}\right)\right).$$
Lemma \ref{convDSBP} yields $\tt\left(L_tf\left(L_{s_1},\ldots,L_{s_p}\right)\right)\leq\tt\left(L_sf\left(L_{s_1},\ldots,L_{s_p}\right)\right)$ since $f$ is bounded and $X^{1/n}_u\leq L_u$ $\tt$ a.s. for every $u\geq0$.\cq

Let us set, for every $t\geq0$,
$$\wt{G}_t=\bigcap_{s>t}G_s.$$  
Recall that $D$ denotes the set of positive dyadic numbers. From Lemma \ref{Lsurmart} and classical results on supermartingales, we can define a right-continuous supermartingale $(\wt{L}_t,t\geq0)$ with respect to the filtration $(\wt{G}_t,t\geq0)$ by setting, for every $t\geq0$,
\be\label{defwtL_t}
\wt{L}_t=\build{\lim}_{s\downarrow t,s\in D}^{}L_s
\ee
where the limit holds $\tt$ a.s. and in $\L^1$. In a way similar to Section \ref{sec3} we can prove that $(\wt{L}_t,t\geq0)$ is a c\`adl\`ag modification of $(L_t,t\geq0)$ with no fixed discontinuities.

From now on, to simplify notation we replace $(L_t,t\geq0)$ by its c\`adl\`ag modification $(\wt{L}_t,t\geq0)$. 

\begin{proposition}
$(L_t,t\geq0)$ is a DSBP which becomes extinct $\tt$ a.s.
\end{proposition}

\proof By the same arguments as in the proof of (\ref{conv2L_t1}), we can prove that, for every $0<s<t$, the following convergence holds in probability under $\tt$,
\be\label{convprobaL_t}
Z\left(\frac{[nt]-[ns]}{n},\frac{[nt]-[ns]+1}{n}\right)\build{\longrightarrow}_{n\to\infty}^{}L_{t-s}.
\ee
Let $s,t,s_1,\ldots,s_p\geq0$ such that $0\leq s_1\leq\ldots\leq s_p\leq s<t$, $\lm>0$ and let $f:\R^p\rightarrow\R$ be a bounded measurable function. For every $n\geq1$, under $\tt^{1/n}$, $(X^{1/n}_k,k\geq0)$ is a Galton-Watson process started at one so that
$$\tt^{1/n}(f(Y^n_{s_1},\ldots,Y^n_{s_p})\exp(-\lm Y^n_{t}))=\tt^{1/n}\left(f(Y^n_{s_1},\ldots,Y^n_{s_p})(\tt^{1/n}(\exp(-\lm X^{1/n}_{[nt]-[ns]})))^{Y^n_s}\right).$$
From Lemma \ref{convDSBP}, (\ref{convprobaL_t}) and dominated convergence, we get
$$\tt(f(L_{s_1},\ldots,L_{s_p})\exp(-\lm L_{t}))=\tt\left(f(L_{s_1},\ldots,L_{s_p})(\tt(\exp(-\lm L_{t-s})))^{L_s}\right).$$ 
Then, $(L_t,t\geq0)$ is a continuous-time Markov chain with values in $\Z_+$ satisfying the branching property. Furthermore, since $\h(\t)<\infty$ $\tt$ a.s., it is immediate that $(L_t,t\geq0)$ becomes extinct $\tt$ a.s. \cq

\subsection{Identification of the probability measure $\tt$}

Let us now define, for every $\t\in\T$ and $t\geq0$, $N_t(\t)=\#\{\sg\in\t:d(\rho,\sg)=t\}$ where we recall that $\rho$ denotes the root of $\t$.

\begin{proposition}\label{nbindividus}
For every $t\geq0$, $N_t=L_t$ $\tt$ a.s.
\end{proposition}
Note that for every $t\geq0$, $L_t$ is the number of subtrees of $\t$ above level $t$. 

\proof Since $\tt(\h(\t)=0)=0$, we have $L_0=1=N_0$ $\tt$ a.s. Thanks to Propositions \ref{defL_t2}, \ref{conv2L_t} and \ref{L_tcadlag}, for every $t>0$, $\tt$ a.s., there exists $h_0>0$ such that for every $h\in(0,h_0]$, $L_t=L_{t-h}=L_{t+h}=Z(t-h,t+h)$. 

The remaining part of the argument is deterministic. We fix $t,h_0>0$ and a (deterministic) tree $\t\in\T$. We assume that there is a positive integer $p$ such that for every $h\in(0,h_0]$,
$$L_t=L_{t-h}=L_{t+h}=Z(t-h,t+h)=p,$$
and we will verify that $N_t(\t)=p$. We denote by $\t^1,\ldots,\t^p$ the $p$ subtrees of $\t$ above level $t-h_0$ and we write $\rho_i$ for the root of the subtree $\t^i$. For every $i\in\{1,\ldots,p\}$, we have $\h(\t^i)>2h_0$ so that there exists $x_i\in\t^i$ such that $d(\rho_i,x_i)=2h_0$. Let us prove that for every $i\in\{1,\ldots,p\}$, 
\be\label{segments}
\t^i_{\leq2h_0}=\{\sg\in\t^i:d(\rho_i,\sg)\leq2h_0\}=[[\rho_i,x_i]].
\ee
To this end, we argue by contradiction and assume that we can find $i\in\{1,\ldots,p\}$ and $t_i\in\t^i$ such that $d(\rho_i,t_i)\leq2h_0$ and $t_i\notin[[\rho_i,x_i]]$. Let $z_i$ be the unique vertex of $\t^i$ satisfying $[[\rho_i,z_i]]=[[\rho_i,x_i]]\cap[[\rho_i,t_i]]$. We choose $c>0$ such that $d(\rho_i,z_i)<c<d(\rho_i,t_i)$. Then it is not difficult to see that $\t$ has at least $p+1$ subtrees above level $t-h_0+c$. This is a contradiction since $L_{t-h_0+c}=p$. So $N_t(\t)=p$, which completes the proof.\cq

Proposition \ref{nbindividus} means that $(L_t,t\geq0)$ is a modification of the process $(N_t,t\geq0)$ which describes the evolution of the number of individuals in the tree. Let us denote by $Q$ the generator of $(L_t,t\geq0)$ which is of the form
$$Q=\left(\begin{array}{cccccc} 0 & 0 & 0 & 0 & 0 & \ldots \\
                               a\gamma(0) & -a & a\gamma(2) & a\gamma(3) & a\gamma(4) & \ldots \\
                               0 & 2a\gamma(0) & -2a & 2a\gamma(2) & a\gamma(3) & \ldots \\
                               0 & 0 & 3a\gamma(0) & -3a& 3a\gamma(2) & \ldots \\
                               \vdots & \vdots & \ddots & \ddots & \ddots & \ddots
                               \end{array}\right),$$
where $a>0$ and $\gamma$ is a critical or subcritical offspring distribution with $\gamma(1)=0$. 

For every $t\geq0$ we let $\f_t$ be the $\sg$-field on $\T$ generated by the mapping $\t\lmt\t_{\leq t}$ and completed with respect to $\tt$. Thus $(\f_t,t\geq0)$ is a filtration on $\T$.
  
\begin{lemma}\label{Mf}
Let $t>0$ and $p\in\N$. Under $\tt$, conditionally on $\f_t$ and given $\{L_t=p\}$, the $p$ subtrees of $\t$ above level $t$ are independent and distributed according to $\tt$. 
\end{lemma}

\proof Thanks to Lemma \ref{discretisation} and Lemma \ref{arbreGWbis}, we can construct on the same probability space $(\OO,\PP)$, a sequence of $\T$-valued random variables $(\TT_n)_{n\geq1}$ distributed according to $\tt^{1/n}$ and a sequence of $\a$-valued random variables $(\ttt_n)_{n\geq1}$ distributed according to $\Pi_{\mu_{1/n}}$ such that, for every $n\geq1$,
\be\label{idcasfini2}
d_{GH}\left(\TT_n,n^{-1}\t^{\ttt_n}\right)\leq4n^{-1}.
\ee
For every $n\geq1$ and $k\geq0$, we define $\XX^n_k=\#\{u\in\ttt_n:|u|=k\}$. Let $t\geq0$ and $p\geq1$, let $g:\T\rightarrow\R$ be a bounded continuous function and let $G:\T^p\rightarrow\R$ be a bounded continuous symmetric function. For $n\geq1$, on the event $\{\XX^n_{[nt]}=p\}$, we set $\{u^n_1,\ldots,u^n_p\}=\{u\in\ttt_n:|u|=[nt]\}$ and $\ttt_n^{i}=\tau_{u^n_i}\ttt_n$ for every $i\in\{1,\ldots,p\}$. Then we can write, thanks to the branching property of Galton-Watson trees,
\begin{eqnarray}
&&\hspace{-1cm}\EE\left(\ind{\left\{\XX^n_{[nt]}=p\right\}}g\left(n^{-1}\t^{\ttt_n}_{\leq [nt]}\right)G\left(n^{-1}\t^{\ttt_n^1},\ldots,n^{-1}\t^{\ttt_n^p}\right)\right)\label{idcasfini1}\\
&=&\EE\left(\ind{\left\{\XX^n_{[nt]}=p\right\}}g\left(n^{-1}\t^{\ttt_n}_{\leq [nt]}\right)\right)(\Pi_{\mu_{1/n}})^{\otimes p}\left(G(n^{-1}\t^{\theta_1},\ldots,n^{-1}\t^{\theta_p})\right),\nonumber
\end{eqnarray}
where $\theta_1,\ldots,\theta_p$ denote the coordinate variables under the product measure $(\Pi_{\mu_{1/n}})^{\otimes p}$. As a consequence of (\ref{idcasfini2}), we see that the law of $n^{-1}\t^\theta$ under $\Pi_{\mu_{1/n}}$ converges to $\tt$ in the sense of weak convergence of measures  on the space $\T$. Then, thanks to Lemma \ref{convDSBP}, the right-hand side of (\ref{idcasfini1}) converges as $n\la\infty$ to 
$$ \tt\left(\ind{\{L_t=p\}}g(\t_{\leq t})\right)\tt^{\otimes p}(G(\t_1,\ldots,\t_p)).$$
Similarly, the left-hand side of (\ref{idcasfini1}) converges as $n\la\infty$ to 
$$ \tt\left(\ind{\{L_t=p\}}g(\t_{\leq t})G(\t^1,\ldots,\t^p)\right),$$
where $\t^1,\ldots,\t^p$ are the $p$ subtrees of $\t$ above level $t$ on the event $\{L_t=p\}$. This completes the proof.\cq

Let us define $J=\inf\{t\geq0:L_t\neq1\}$. Then $J$ is an $(\f_t)_{t\geq0}$-stopping time. 

\begin{lemma}\label{MF}
Let $p\in\N$. Under $\tt$, given $\{L_J=p\}$, the $p$ subtrees of $\t$ above level $J$ are independent and distributed according to $\tt$, and are independent of $J$.
\end{lemma}

\proof Let $p\in\N$, let $f:\R_+\rightarrow\R$ be a bounded continuous function and let $G:\T^p\rightarrow\R$ be a bounded continuous symmetric function. On the event $\{L_J=p\}$, we denote by $\t^1,\ldots,\t^p$ the $p$ subtrees of $\t$ above level $J$. Let $n\geq1$ and $k\geq0$. On the event $\{L_{(k+1)/n}=p\}$, we denote by $\t^{1,(n,k)},\ldots,\t^{p,(n,k)}$ the $p$ subtrees of $\t$ above level $(k+1)/n$. On the one hand, the right-continuity of the mapping $t\lmt L_t$ gives
\ba
&&\hspace{-3.5cm}\tt\left(\sum_{k=1}^\infty\ind{\left\{L_{(k+1)/n}=p\right\}}G(\t^{1,(n,k)},\ldots,\t^{p,(n,k)})f\left((k+1)/n\right)\ind{\left\{k/n<J\leq(k+1)/n\right\}}\right)\\
&\build{\la}_{n\to\infty}^{}&\tt\left(\ind{\{L_J=p\}}G(\t^1,\ldots,\t^p)f(J)\right).\ea
On the other hand, thanks to Lemma \ref{Mf}, we can write for every $n\geq1$ and $k\geq0$,
\ba
&&\hspace{-1cm}\tt\left(\ind{\left\{L_{(k+1)/n}=p\right\}}G(\t^{1,(n,k)},\ldots,\t^{p,(n,k)})f\left((k+1)/n\right)\ind{\left\{k/n<J\leq(k+1)/n\right\}}\right)\\
&=&\tt\left(\ind{\left\{L_{(k+1)/n}=p\right\}}f\left((k+1)/n\right)\ind{\left\{k/n<J\leq(k+1)/n\right\}}\right)\tt^{\otimes p}(G(\t_1,\ldots,\t_p)).
\ea
It follows that $\tt\left(\ind{\{L_J=p\}}G(\t^1,\ldots,\t^p)f(J)\right)=\tt\left(\ind{\{L_J=p\}}f(J)\right)\tt^{\otimes p}(G(\t_1,\ldots,\t_p))$.\cq

We can now complete the proof of Theorem \ref{casfini}. The random variable $J$ is the first jump time of the DSBP $(L_t,t\geq0)$ so that $J$ is distributed according to the exponential distribution with parameter $a$ and is independent of $L_J$. Thanks to Proposition \ref{nbindividus}, there exists $\sg_J\in\t$ such that $\t_{\leq J}=[[\rho,\sg_J]]$. Lemma \ref{MF} gives the last part of the description of $\tt$.

Another way to describe $\tt$ is as follows: Assume that we are given on the same probability space $(\OO,\PP)$ an $\a$-valued random variable $\ttt$ distributed according to $\Pi_\gamma$ and an independent sequence of independent random variables $({\bf h}_u,u\in U)$ with values in $[0,\infty)$, such that each variable ${\bf h}_u$ is distributed according to the exponential distribution with parameter $a$. We set ${\bf T}=(\ttt,\{{\bf h}_u\}_{u\in\ttt})$ and $\TT=\t^{\bf{T}}$. Then the random variable $\TT$ is distributed according to $\tt$.


\begin{thebibliography}{99}

\bibitem{Al1} {\sc Aldous, D.} (1991) The continuum random tree I. {\it Ann. Probab.} 
 {\bf 19}, 1-28.

\bibitem{Al3} {\sc Aldous, D.} (1993) The continuum random tree III. {\it Ann. Probab.} 
 {\bf 21}, 248-289.

\bibitem{AN} {\sc Athreya, K.B., Ney, P.E.} (1972) {\it Branching Processes.} Springer, Berlin.

\bibitem{BBI} {\sc Burago, D., Burago, Y., Ivanov, S.} (2001) {\it A Course in Metric Geometry.} Graduate studies in mathematics, vol.33, AMS, Boston.

\bibitem{DM} {\sc Dellacherie, C., Meyer, P.A.} (1980) {\it Probabilit\'es et Potentiels, Chapitres V \`a VIII: Th\'eorie des Martingales.} Hermann, Paris.

\bibitem{DuLG-Mon} {\sc Duquesne, T., Le Gall, J.F.} (2002) Random Trees, L\'evy Processes and Spatial Branching Processes, {\it Ast\'erisque} {\bf 281}.  

\bibitem{DuLG} {\sc Duquesne, T., Le Gall, J.F.} (2005) Probabilistic 
and fractal aspects of L\'evy trees.
{\it Probab. Th. Rel. Fields}, {\bf 131}, 553-603.

\bibitem{EPW} {\sc Evans, S.N., Pitman, J.W., Winter, A.} (2003) 
Rayleigh processes, real trees
and root growth with re-grafting. {\it Probab. Th. Rel. Fields}, to appear.

\bibitem{EW} {\sc Evans, S.N., Winter, A.} Subtree prune and re-graft: A reversible tree valued Markov process. {\it Ann. Probab.}, to appear.

\bibitem{Gr} {\sc Grimvall, A.} (1974) On the convergence of sequences of branching processes. {\it Ann. Probab.} {\bf 2} 1027-1045.

\bibitem{La} {\sc Lamperti, J.} (1966) The limit of a sequence of branching processes. {\it  Z. Wahrsch. verw. Gebiete} {\bf 7} 271-288.

\bibitem{M1} {\sc Miermont, G.} (2003) Self-similar fragmentations derived from the stable tree I: splitting at heights. {\it Probab. Th. Rel. Fields}, {\bf 127}, 423-454.

\bibitem{M2} {\sc Miermont, G.} (2004) Self-similar fragmentations derived from the stable tree II: splitting at nodes. {\it Probab. Th. Rel. Fields}, to appear.
 
\bibitem{No} {\sc Norris, J.} (1997) {\it Markov Chains.} 
Cambridge University Press, Cambridge.

\end{thebibliography}
\end{document}